\documentclass[a4,12pt]{amsart}
\usepackage{amsxtra,amsmath,amsthm,amssymb,amscd,amsfonts,
multicol,tabularx,latexsym}
\theoremstyle{definition}
\newtheorem{propI}{Proposition I}
\newtheorem{propII}{Proposition II}

\newtheorem{thma}{Theorem A}

\theoremstyle{definition}
\newtheorem{thm}{Theorem}
\newtheorem{lem}{Lemma}
\newtheorem{prop}{Proposition}

\newcommand{\Gal}{\mathrm{Gal}}

\newcommand{\Hom}{\mathrm{Hom}}

\newcommand{\rank}{\mathrm{rank}}

\newcommand{\Q}{\mathbb Q}
\newcommand{\Z}{\mathbb Z}
\newcommand{\N}{\mathbb N}
\newcommand{\F}{\mathbb F}

\newcommand{\coker}{\mathrm{coker}}
\newcommand{\im}{\mathrm{im}}

\newcommand{\m}{\mu}
\renewcommand{\t}{\widetilde}

\newcommand{\ab}{\mathrm{ab}}
\newcommand{\cl}{\mathrm{cl}}
\newcommand{\Cl}{\mathrm{Cl}}
\date{}
\begin{document}
\title[]
{Construction of maximal unramified $p$-extensions with
prescribed Galois groups}
\thanks{
\\
This research is partially supported by the Grant-in-Aid for 
Young Scientists (B), Ministry of Education, Science, Sports and Culture,
Japan.}
\author{Manabu\ OZAKI}
\maketitle
\section{Introduction}
For any number field $F$ (not necessary of finite degree)
and prime number $p$, let $L_p(F)$
denote the maximal unramified $p$-extension over $F$, and put
$\t{G}_F(p)=\Gal(L_p(F)/F)$.
Though the structure of $\t{G}_F(p)$ has been
one of the most fascinating theme of number theory,
our knowledge on it is not enough even at present:
It had been a cerebrated open problem for a long time
whether $\t{G}_F(p)$ can be infinite for a number field $F$
of finite degree,
and Golod and Shafarevich solved it b
y giving
$F$ with infinite $\t{G}_F(p)$.
Hence we barely know that $\t{G}_F(p)$ can be infinite.
However we do not know exactly what kind of pro-$p$-groups occur
as $\t{G}_F(p)$; for example, there are no examples of infinite $\t{G}_F(p)$
for number fields $F$ of finite degree
whose structure is completely determined.
The known general property of the group $\t{G}_F(p)$
for number fields $F$ of finite degree 
is only that $\t{G}_F(p)$ is a finitely presented pro-$p$-group 
any whose open subgroup has finite abelianization,
which comes from
rather fundamental facts of algebraic number theory, namely,
class field theory and the finiteness of the
ideal class group.
Indeed, a consequence of the Fontaine-Mazur conjecture
predicts that $\t{G}_F(p)$ has a certain distinguished property
(see \cite{FM}),
however, this conjecture seems far reach object at present.
\par
On the other hand,
we have known that various kind of pro-$p$-groups, especially,
finite $p$-groups in fact occur as $\t{G}_F(p)$.
For example, Scholz and Taussky \cite{ST} have already determined
the structure of $\t{G}_F(p)$ for $F=\Q(\sqrt{-4027})$ and $p=3$
in 1930's :
this group is a non-abelian finite group of order $3^5$.
Also, Yahagi \cite{Yah78} showed
that for any given finite abelian $p$-group $A$
there exists an number field $F$ of finite degree
such that $\t{G}_F(p)^{\mathrm{ab}}\simeq A$.
\par
In the present paper, we shall show that
every {\it finite} $p$-group occurs as $\t{G}_F(p)$ for some
number field $F$ of finite degree:
\begin{thm}\label{thm1}
Let $p$ be a prime number.
Then for any given finite $p$-group $G$,\ there exists a number field
$F$ of finite degree such that $\t{G}_F(p)\simeq G$.
\end{thm}
In the case where we take account number fields of infinite degree,
we shall determine completely the set of the isomorphism classes
of the pro-$p$-groups $\t{G}_F(p)$ where $F$ is an algebraic extension over
$\Q$:
\begin{thm}\label{thm2}
For a prime number $p$, let $\mathcal{C}_p$
be the set of all the isomorphism classes of the pro-$p$-groups
$\t{G}_F(p)$ where $F$ is an algebraic extension over
$\Q$.
Then $\mathcal{C}_p$ is exactly equal to the set
of all the isomorphism classes of
the pro-$p$-groups with countably many generators.
\end{thm}
%
%
\section{Strategy of construction}
We first introduce some notations.
For any number field $F$ and prime number $p$,
let $L_p(F)/F$ and $L^{\ab}_p(F)/F$
be the maximal unramified $p$-extension and the maximal
unramified abelian $p$-extension of $F$, respectively.
Also we write $\Cl(F)$ for the ideal class group of $F$.
For any pro-$p$-group $G$ we define $d(G)$ and $r(G)$
to be the generator rank and relation rank of $G$, respectively;
$d(G)=\dim_{\F_p}H^1(G,\F_p),\ r(G)=\dim_{\F_p}H^2(G,\F_p)$. 
If $\t{G}_F(p)$ is finite, we put 
\begin{equation*}
\begin{aligned}
&B_p(F)\!=\!\left((2[F(\mu_p):F]+1)\left(\#\t{G}_F(p)-1\right)
+2\right)\times\\
&\left(
d(\t{G}_F(p))\!+\!r(\t{G}_F(p))
\!+\!\dim_{\F_p}(\Cl\left(L_p(F)(\mu_p))\otimes\F_p\right)
_{\Gal(L_p(F)(\mu_p)/F(\mu_p))}\right)
\\
&+d(\t{G}_F(p))+4\dim_{\F_p}(\Cl\left(L_p(F)(\mu_p))\otimes\F_p\right)
_{\Gal(L_p(F)(\mu_p)/F(\mu_p))}+3,
\end{aligned}
\end{equation*}
$\mu_p$ being the group of the $p$-th roots of unity,
which is a constant depending only on $F$ and $p$.
In what follows we shall show the following two propositions:
\begin{propI}
Let $p$ be a prime number.
Then there exists a totally imaginary number field $k$
of finite degree such that $\t{G}_{k}(p)=1$
and that $r_2(k)\ge B_p(k)$, where $r_2(k)$ stands for the number of
complex archimedean places of $k$.
\end{propI}
\begin{propII}
Let $p$ be a prime number and
$k$ a totally imaginary number field of finite degree
with finite $\t{G}_k(p)$
such that
$r_2(k)\ge B_p(k)$ and $L_p(k)/k$ is $p$-decomposed, namely,
every prime of $k$ lying over $p$
splits completely in $L_p(k)$.
Then for any exact sequence of  finite $p$-groups
\[
1\longrightarrow \Z/p\longrightarrow G'\overset{\pi}{\longrightarrow}
\t{G}_k(p)\longrightarrow 1
\]
there exists a finite extension $k'/k$
such that $L_p(k)\cap k'=k$,
$r_2(k')\ge B_p(k')$, $L_p(k')/k'$ is $p$-decomposed,
and $\t{G}_{k'}(p)$ fits the following commutative diagram
with an isomorphism
$\theta:G'\overset{\sim}{\rightarrow}\t{G}_{k'}(p)$:
\begin{equation}\label{cd-propII}
\begin{CD}
 G' @>\pi>> \t{G}_k(p)\\
@V{\theta\ \wr}VV @| \\
\t{G}_{k'}(p) @>{\mathrm{restriction}}>> \t{G}_k(p),
\end{CD}
\end{equation}
where the bottom horizontal map is the restriction map
induced by the inclusion $L_p(k)\subseteq L_p(k')$.
\end{propII}
Then one can derive Theorems 1 and 2 from Propositions I and II as
follows:
Theorem 1 follows from Propositions I and II by induction on the
order of the finite $p$-group $G$ because
any non-trivial finite $p$-group has a normal subgroup of order $p$.
Theorem \ref{thm2} also follows similarly from Propositions I and II:
Because the Galois group of a pro-$p$-extension of number fields
should have countably many generators,
it is enough to show that for any given pro-$p$ group $G$
with countably many generators there exists number field $F$
with $\t{G}_F(p)\simeq G$.
We may assume that $G$ is infinite thanks to Theorem 1.
Then there exists a projective system of finite $p$-groups 
$G_n\ (0\le n\in\Z)$ with $G_0=1$ and surjective morphisms $\phi_{m,n}:
G_m\rightarrow G_n\ (m\ge n)$ with $\ker \phi_{n+1,n}\simeq\Z/p$ such that 
$G$ is isomorphic to the projective limit of
this projective system.
Hence, by using Propositions I and II,
we can construct a tower of number fields $k_n$ of finite degree
\[
k_0\subseteq k_1\subseteq\cdots\subseteq k_n\subseteq k_{n+1}\subseteq\cdots
\subseteq F:=\bigcup_{n\in\N}k_n
\]
such that $L_p(k_n)\cap F=k_n$
and $\t{G}_{k_n}(p)$'s fit the commutative diagram
\begin{equation*}
\begin{CD}
G_{n+1} @>{\phi_{n+1,n}}>> G_n\\
@V{\theta_{n+1}\ \wr}VV @V{\theta_n\ \wr}VV\\
\t{G}_{k_{n+1}}(p) @>{\mathrm{restriction}}>> \t{G}_{k_n}(p)
\end{CD}
\end{equation*}
with suitable isomorphisms $\theta_n$'s.
Thus $F$ is a required field because $L_p(F)=\bigcup_{n\in\N} L_p(k_n)$
and $\t{G}_F(p)\simeq\varprojlim \t{G}_{k_n}(p)\simeq G$.

We can easily show Proposition I thanks to Horie's theorem \cite{Hor}
on the indivisibility of the class numbers of imaginary quadratic fields
and the Ferrero-Washington theorem \cite{FW} on
the vanishing of the Iwasawa $\mu$-invariants of abelian number fields.
However our proof of Proposition II is rather complicated and long,
in which we shall employ extensively the Chebotarev density theorem
and the theory of embedding problem of Galois extensions.
%
%
%
%
%
%
%
%
%
\section{Proof of Proposition I}
%
%
%
%
%
In this section, we shall give a proof of Proposition I.

Let $p$ be any prime number.
Then it follows Horie \cite{Hor} that
there exists an imaginary quadratic field $F$
such that the prime number $p$ does not decompose
in $F$ and does not divide the class number $h_F$ of $F$.
Then the class number of the $n$-th layer $F_n$ of the cyclotomic
$\Z_p$-extension $F_\infty/F$ is prime to $p$ by Iwasawa \cite{Iw56a}
as well known, since $F_\infty$ has the unique prime
lying over $p$ by $p\nmid h_F$.
Hence we have $\t{G}_{F_n}(p)=1$, from which we see that 
\[
B_p(F_n)=6\dim_{\F_p}\left(\Cl(F_n(\mu_p))\otimes\F_p\right)+3.
\]
On the other hand, $\dim_{\F_p}(\Cl\left(F_n(\mu_p))\otimes\F_p\right)$
is bounded as $n$ goes to infinity because
the Iwasawa $\mu$-invariant of the cyclotomic $\Z_p$-extension
$F_\infty(\mu_p)/F(\mu_p)$,
whose $n$-th layer is $F_n(\mu_p)$,
vanishes by the Ferrero-Washington theorem \cite{FW}.
Therefore inequality $r_2(F_n)\ge B_p(F_n)$ holds
for sufficiently large $n$, which proves Proposition I.
%
%
%
%
%
\section{Preliminary results on the structure
of the group of units}
Since our proof of Proposition II is rather long and complicated,
we divide it into 5 sections.
In this section we shall give preliminary results on
the Galois module structure of the group of units modulo $p$-th powers,
which we shall need below.

In what follows, we fix a prime number $p$ and totally imaginary number field
$k$ of finite degree with finite $\t{G}_k(p)$.
We first introduce some notations.
For any number field $F$ we write $\mathcal{O}_F$
for the integer ring of $F$. 
We put $G:=\t{G}_k(p)=\Gal(L_p(k)/k)$,
$\Delta:=\Gal(k(\mu_p)/k)$ and $K:=L_p(k)(\mu_p)$.
Then the order of $\Delta$ divides $p-1$,
and we have the natural isomorphism
\[
\Gal(K/k)\simeq\Delta\times G
\]
induced by the restriction of field automorphisms and 
we shall identify these two groups via the above isomorphism. 
For any $\Z_p[\Delta]$-module $M$ and a $p$-adic character 
$\chi\in\widehat{\Delta}:=\Hom(\Delta,\Z_p^\times)$,
we put 
\[
e_\chi:=\frac{1}{\#\Delta}
\sum_{\delta\in\Delta}\chi(\delta)\delta^{-1}
\in\Z_p[\Delta],
\]
and define the $\chi$-part $M^\chi$ of $M$
to be $e_\chi M$.
Then we have the direct sum decomposition
\[
M=\bigoplus_{\chi\in\widehat{\Delta}}M^\chi
\]
as a $\Z_p[\Delta]$-module.
We define the cyclotomic character $\omega\in\widehat{\Delta}$
by $\delta(\zeta)
=\zeta^{\omega(\delta)}\ (\delta\in\Delta,\ \zeta\in\mu_p)$.
\begin{lem}\label{chi-unit}
For $\chi\in\widehat{\Delta}$, we have
\begin{equation*}
\dim_{\F_p}(\mathcal{O}_K^\times\otimes\F_p)^\chi
=
\begin{cases}
 r_2(k)\# G\ \ (\mbox{if $\chi\ne 1,\omega$ or $\Delta=1$}),\\
 r_2(k)\# G+1\ (\mbox{if $\chi=\omega\ne 1$}),\\
 r_2(k)\# G-1\ (\mbox{if $\chi=1\ne\omega$}).
\end{cases}
\end{equation*}
\end{lem}
\noindent
{\bf Proof.}\ \
It follows from Herbrand's unit theorem
(see \cite[Chapter V, Proposition 2.3]{Jan} for example)
that
\begin{equation*}
 \mathcal{O}_K^\times\otimes\Q_p\simeq
\Q_p[\Delta\times G]^{\oplus r_2(k)-1}\oplus
I_{\Q_p[\Delta\times G]}
\end{equation*}
as $\Gal(K/k)(=\Delta\times G)$-modules,
where $I_R$ denotes the augmentation ideal of $R$
for any group ring $R$.
Hence we have
\begin{equation*}
(\mathcal{O}_K^\times\otimes\Q_p)^\chi\simeq
\begin{cases}
\Q_p[G]^{\oplus r_2(k)-1}\oplus
I_{\Q_p[G]}\ (\mbox{if $\chi=1$}),\\
\Q_p[G]^{\oplus r_2(k)}\ (\mbox{if $\chi\ne 1$}),
\end{cases}
\end{equation*}
as $G$-modules, which implies
\begin{equation*}
\rank_{\Z_p}(\mathcal{O}_K^\times\otimes\Z_p)^\chi
=
\begin{cases}
r_2(k)\#G-1\ (\mbox{if $\chi=1$}),\\
r_2(k)\#G\ (\mbox{if $\chi\ne 1$}).
\end{cases}
\end{equation*}
Therefore we obtain the lemma by noting
\begin{equation*}
\dim_{\F_p}(\mathcal{O}_K^\times\otimes\F_p)^\chi
=
\begin{cases}
\rank_{\Z_p}(\mathcal{O}_K^\times\otimes\Z_p)^\chi+1\ 
(\mbox{if $\chi=\omega$}),\\
\rank_{\Z_p}(\mathcal{O}_K^\times\otimes\Z_p)^\chi\ 
(\mbox{otherwise}).
\end{cases}
\end{equation*} 
\hfill$\Box$

\

We shall recall some basic facts on the category of the finitely generated
$\F_p[G]$-modules for a finite $p$-group $G$.
For any $\F_p[G]$-module $A$, we put $A^*=\mathrm{Hom}_{\F_p}(A,\F_p)$
and define left $G$-action on it by $\sigma f:=f\circ\sigma^{-1}$
for $f\in A^{*}$ and $\sigma\in G$.
Then we find that $(A^*)^*\simeq A$ 
and $\F_p[G]^*\simeq\F_p[G]$ as $\F_p[G]$-modules
for any finitely generated $\F_p[G]$-module $A$,
from which we see that every finitely generated free $\F_p[G]$-module
is an injective object in the category of finitely generated
$\F_p[G]$-modules.
 Therefore every free $\F_p[G]$-submodule of a finitely generated
$\F_p[G]$-module is a direct factor of it.
\par
For any finitely generated $\F_p[G]$-module $A$, there exists a surjection
\[
\F_p[G]^{\oplus r}\longrightarrow A
\]
as $\F_p[G]$-modules
for $r=\dim_{\F_p}A_G$ by Nakayama's lemma.
Replacing $A$ by $A^*$ and taking the dual $*$ in the above, we find that there exists an injection
\[
A\longrightarrow \F_p[G]^{\oplus r^*}
\]
as $\F_p[G]$-modules
for $r^*=\dim_{\F_p}A^G$ since $(A^*)_G\simeq(A^G)^*\simeq A^G$
as $\F_p$-modules.
\begin{lem}\label{lem2}
For any $\chi\in\widehat{\Delta}$,
there exists a direct sum decomposition 
\[
(\mathcal{O}_K^\times\otimes\F_p)^\chi=F_\chi\oplus N_\chi 
\]
as $\F_p[G]$-modules
such that
\begin{equation*}
\begin{aligned}
&F_\chi\simeq\F_p[G]^{\oplus r_\chi},\ N_GN_\chi=0,\\
&r_\chi\ge r_2(k)-(\#G-1)(d(G)+r(G)+\dim_{\F_p}(\Cl(K)\otimes\F_p)_G)
-\delta_{G,\chi},\\
&\!\dim_{\F_p}\!N_\chi^G\!\le\!
([k(\mu_p)\!:k](\# G\!-\!1)\!+\!1)(d(G)\!+\!r(G)\!+\!\dim_{\F_p}\!
(\Cl(K)\!\otimes\!\F_p)_G),
\end{aligned}
\end{equation*}
where $N_G=\sum_{g\in G}g\in\F_p[G]$,
and $\delta_{G,\chi}=1$ or 0 according to [$G=1$ and $\chi=1$] or not. 
\end{lem}
\noindent{\bf Proof.}\ \ \ 
Let $F_\chi$ be a maximal free $\F_p[G]$-submodule of
$(\mathcal{O}_K^\times\otimes\F_p)^\chi$.
Then there exists a direct sum decomposition
\begin{equation}\label{dsum}
(\mathcal{O}_K^\times\otimes\F_p)^\chi
=F_\chi\oplus N_\chi
\end{equation}
by the injectivity of the free $\F_p[G]$-modules.
Since $N_\chi$ has no non-trivial free $\F_p[G]$-submodule,
the annihilator $\mathrm{Ann}_{\F_p[G]}(x)$ of each $x\in N_\chi$
is a non-zero left ideal of $\F_p[G]$.
Hence we have $\F_pN_G=\F_p[G]^G\supseteq
(\mathrm{Ann}_{\F_p[G]}(x))^G\ne 0$ for any $x\in N_\chi$,
which implies $N_G\in\mathrm{Ann}_{\F_p[G]}(x)$ and $N_GN_\chi=0$.

Finally, we shall estimate the free rank $r_\chi$ of $F_\chi$
and $\dim_{\F_p}N_\chi^G$.
We derive from the exact sequence
\[
0\longrightarrow\mathcal{O}_K^\times/\mu_p
\overset{p}{\longrightarrow}
\mathcal{O}_K^\times
\longrightarrow\mathcal{O}_K^\times\otimes\F_p\longrightarrow 0
\]
the $G$-cohomology exact sequence (identifying $G$ with $\Gal(K/k(\mu_p)$))
\begin{equation}\label{seq1}
(\mathcal{O}_K^\times/\mu_p)^G
\overset{p}{\longrightarrow}
\mathcal{O}_{k(\mu_p)}^\times
\longrightarrow(\mathcal{O}_K^\times\otimes\F_p)^G
\longrightarrow H^1(G,\mathcal{O}_K^\times/\mu_p).
\end{equation}
On the other hand,
from the exact sequence
\[
0\longrightarrow\mu_p
\longrightarrow\mathcal{O}_K^\times
\longrightarrow\mathcal{O}_K^\times/\mu_p
\longrightarrow 0,
\]
we get the exact sequence
\begin{equation}\label{seq2}
H^1(G,\mathcal{O}_K^\times)
\longrightarrow
H^1(G,\mathcal{O}_K^\times/\mu_p)
\longrightarrow
H^2(G,\mu_p)=H^2(G,\F_p).
\end{equation}
Then it follows from \eqref{seq1} and \eqref{seq2} that
\begin{equation}\label{est}
\begin{split}
\dim_{\F_p}(\mathcal{O}_K^\times\otimes\F_p)^G
\leq r_2(k)[k(\mu_p):k]+\dim_{\F_p}
(H^1(G,\mathcal{O}_K^\times/\mu_p)\otimes\F_p)
\\\leq r_2(k)[k(\mu_p):k]
+\dim_{\F_p}(H^1(G,\mathcal{O}_K^\times)\otimes\F_p)
+r(G).
\end{split}
\end{equation}

If we let $j:\mathrm{Cl}(k(\mu_p))\longrightarrow\mathrm{Cl}(K)$
be the natural map from the ideal class group of $k(\mu_p)$
to that of $K$ induced by the inclusion $k(\mu_p)\subseteq K$,
then 
\begin{equation}\label{H^1}
H^1(G,\mathcal{O}_K^\times)\simeq\mathrm{ker}\,j
\subseteq\Cl(k(\mu_p))\otimes\Z_p
\end{equation}
as well known (see \cite{Iw56b} for example).

Let $L$ be the maximal unramified
elementary abelian $p$-extension
of $K$ such that $L/k(\m_p)$ is a Galois extension
and $\Gal(K/k(\mu_p))$ acts on $\Gal(L/K)$
trivially via the inner automorphisms of $\Gal(L/k(\mu_p))$,
and let $L_0$ be the maximal unramified elementary abelian
$p$-extension of $k(\mu_p)$,
then $\Gal(L/K)\simeq(\Cl(K)\otimes\F_p)_G$
and we have
\begin{align*}
\dim_{\F_p}(\Cl(k(\mu_p))\otimes\F_p)
&=
\dim_{\F_p}\Gal(L_0/k(\mu_p))\\
&\le
\dim_{\F_p}\Gal(K/k(\mu_p))\otimes\F_p
+
\dim_{\F_p}\Gal(L/K)\\
&=
d(G)+\dim_{\F_p}(\Cl(K)\otimes\F_p)_G,
\end{align*}
since $L_0K\subseteq L$.
Hence we obtain
\begin{equation}\label{^G}
\dim_{\F_p}(\mathcal{O}_K^\times\otimes\F_p)^G
\leq
r_2(k)[k(\mu_p):k]
+d(G)+r(G)+\dim_{\F_p}(\Cl(K)\otimes\F_p)_G
\end{equation}
from \eqref{est} and \eqref{H^1}.
Also we have
\begin{equation}\label{ineq1}
\begin{aligned}
\dim_{\F_p}(\mathcal{O}_K^\times\otimes\F_p)^\chi
&=r_\chi\#G+\dim_{\F_p}N_\chi\\
&=r_\chi\#G+\dim_{\F_p}N_\chi^*
\end{aligned}
\end{equation}
by \eqref{dsum}.

Let $\{n_1,....,n_s\}\subseteq N_\chi^*\ (s=\dim_{\F_p}(N_\chi^*)_G)$
be a set of generators of $N_\chi^*$ over $\F_p[G]$.
Then $\dim_{\F_p}\F_p[G]n_i\leq \#G-1$
since $\mathrm{Ann}_{\F_p[G]}(n_i)\ne 0$.
Hence we obtain 
\begin{equation}\label{ineq2}
\begin{aligned}
\dim_{\F_p}N_\chi^*&\le
\sum_{i=1}^s\dim_{\F_p}\F_p[G]n_i\le(\#G-1)\dim_{\F_p}(N_\chi^*)_G\\
=&(\#G-1)\dim_{\F_p}N_\chi^G
=(\#G-1)(\dim_{\F_p}((\mathcal{O}_K^\times\otimes\F_p)^\chi)^G-r_\chi).
\end{aligned}
\end{equation}
Then it follows from \eqref{ineq1} and \eqref{ineq2} that
\begin{equation}\label{ineq3}
\dim_{\F_p}(\mathcal{O}_K^\times\otimes\F_p)^\chi
\le
r_\chi+(\#G-1)\dim_{\F_p}((\mathcal{O}_K^\times\otimes\F_p)^\chi)^G.
\end{equation}

In what follows, we may assume that $G\ne 1$ because
the lemma is obvious in this case by Lemma \ref{chi-unit}.
Suppose that
\begin{equation}\label{contra}
r_\psi<r_2(k)-(\#G-1)(d(G)+r(G)+\dim_{\F_p}(\Cl(K))\otimes\F_p)_G)
\end{equation}
for some $\psi\in\widehat{\Delta}$.
Then by summing up the both side of \eqref{ineq3} over all the 
$\chi\in\widehat{\Delta}$, and by using \eqref{^G} and assumption
\eqref{contra},
we have
\begin{equation*}
\begin{aligned}
&r_2(k)[k(\mu_p):k]\#G
=\dim_{\F_p}(\mathcal{O}_K^\times\otimes\F_p)\\
&\le\sum_{\chi\ne \psi}r_\chi+r_\psi+(\#G-1)
(\dim_{\F_p}(\mathcal{O}_K^\times\otimes\F_p)^G)\\
&<\sum_{\chi\ne \psi}r_\chi+
r_2(k)-(\#G-1)(d(G)+r(G)+\dim_{\F_p}(\Cl(K)\otimes\F_p)_G)\\
&+(\#G-1)(r_2(k)[k(\mu_p):k]+d(G)+r(G)
+\dim_{\F_p}(\Cl(K)\otimes\F_p)_G
)\\
&=
\sum_{\chi\ne \psi}r_\chi+r_2(k)+(\#G-1)r_2(k)[k(\mu_p):k],
\end{aligned}
\end{equation*}
from which we derive
\begin{equation}\label{contr2}
r_2(k)(\#\Delta-1)<\sum_{\chi\ne \psi}r_\chi.
\end{equation}

On the other hand, it follows from Lemma \ref{chi-unit}
and assumption $G\ne 1$ that
\[
r_\chi\le r_2(k)
\]
for any $\chi\in\widehat{\Delta}$.
Hence it follows from \eqref{contr2} that
\[
r_2(k)(\#\Delta-1)<\sum_{\chi\ne \psi}r_\chi
\le r_2(k)(\#\Delta-1),
\]
which is contradiction. Therefore \eqref{contra} does not hold for any
$\psi\in\widehat{\Delta}$, which implies
the lower estimation on $r_\chi$. 

We obtain the following estimation on $\dim_{\F_p}N_\chi^G$
from that on $r_\chi$ and \eqref{^G}:
\begin{equation*}
\begin{aligned}
&\dim_{\F_p}N_\chi^G\le
\sum_{\psi\in\widehat{\Delta}}\dim_{\F_p}N_\psi^G
=
\dim_{\F_p}(\mathcal{O}_K^\times\otimes\F_p)^G
-\sum_{\psi\in\widehat{\Delta}}r_\psi\\
&\le
r_2(k)[k(\mu_p):k]
+d(G)+r(G)+\dim_{\F_p}(\Cl(K)\otimes\F_p)_G\\
&-\#\Delta(r_2(k)-(\#G-1)(d(G)+r(G)+\dim_{\F_p}(\Cl(K)\otimes\F_p)_G))\\
&=
([k(\mu_p):k](\#G-1)+1)(d(G)+r(G)+\dim_{\F_p}(\Cl(K)\otimes\F_p)_G).
\end{aligned}
\end{equation*}
This completes the proof of the lemma.
\hfill$\Box$
%
%
%
%
\section{Application of the Chebotarev density theorem}
%
%
In this section we shall give fundamental of our construction,
which is an application of Chebotarev density theorem.

Under the same setting and notations as in the preceding section,
we shall give the following:
\begin{lem}\label{lem3}
Let 
\[
f:\mathcal{O}_K^\times\otimes\F_p\longrightarrow\F_p[\Delta\times G]^{\oplus r}
\]
be any given $\F_p[\Delta\times G]$-homomorphism
, $M/K$ any finite abelian extension
linearly disjoint from $K(\sqrt[p]{\mathcal{O}_K^\times})/K$,
and $\sigma_1,\dots,\sigma_r\in\Gal(M/K)$ any automorphisms.
Then there exist distinct primes
$\frak{L}_1,\dots,\frak{L}_r$ of $K$ and $g_1,\dots,g_r\in\Z$
with the following properties:\\
(i)\ \ The primes $\frak{L}_1,\dots,\frak{L}_r$ are unramified in $K/\Q$,
of degree one, and lying over distinct rational primes different from $p$.\\
(ii) We have $(\frak{L_i},M/K)=\sigma_i$ for $1\le i\le r$,
where $(*,M/K)$ denotes the Artin symbol for $M/K$.\\
(iii)\ \ $(g_i\mod\frak{L}_i)\otimes 1$ generates
$(\mathcal{O}_K/\frak{L}_i)^\times\otimes\F_p\simeq\F_p$
for each $i$,
and if we define $\F_p[\Delta\times G]$-isomorphism
\[
\varphi:\F_p[\Delta\times G]^{\oplus r}
\longrightarrow
\bigoplus_{i=1}^r\bigoplus_{\sigma\in \Delta\times G}
(\mathcal{O}_K/\sigma\frak{L}_i)^\times\otimes\F_p
\]
by
\[
\varphi\left(
\left(\sum_{\sigma\in \Delta\times G}a_\sigma^{(i)}\sigma\right)_{1\le i\le r}
\right)
=\left(
((g_i\!\!\!\!\!\mod\sigma\frak{L}_i)\otimes a_\sigma^{(i)})_{\sigma\in \Delta\times G}
\right)_{1\le i\le r},
\]
and 
\[
\pi:\mathcal{O}_K^\times\otimes\F_p\longrightarrow
\bigoplus_{i=1}^r\bigoplus_{\sigma\in \Delta\times G}
(\mathcal{O}_K/\sigma\frak{L}_i)^\times\otimes\F_p
\]
to be the natural map induced by the projections 
$
\mathcal{O}_K^\times\longrightarrow
(\mathcal{O}_K/\sigma\frak{L}_i)^\times
$,
then we have
\[
\varphi\circ f=\pi.
\]
\end{lem}
\noindent
{\bf Proof.}\ \ \ 
Let $\mathrm{Pr}_i:\F_p[\Delta\times G]^{\oplus r}
\longrightarrow
\F_p[\Delta\times G]
$
be the $i$-th projection for $1\le i\le r$,
and we fix an $\F_p$-isomorphism $\iota:\F_p\simeq\mu_p$.
Define homomorphism of $\F_p$-modules
\[
\phi_i:
\mathcal{O}_K^\times\otimes\F_p\longrightarrow
\mu_p
\]
by $\phi_i(x)=\iota(a_1^{(i)})$
when $\mathrm{Pr}_i\circ f(x)
=\sum_{\sigma\in \Delta\times G}a_\sigma^{(i)}\sigma$.
Then it follows from the Kummer duality
\[
\Hom(\mathcal{O}_K^\times\otimes\F_p,\mu_p)\simeq
\Gal(K(\sqrt[p]{\mathcal{O}_K^\times})/K)
\]
that there exists the unique element
$\gamma_i\in\Gal(K(\sqrt[p]{\mathcal{O}_K^\times})/K)$
such that
\begin{equation}\label{gamma}
\phi_i(\varepsilon\otimes 1)
=\gamma_i(\sqrt[p]{\varepsilon})(\sqrt[p]{\varepsilon})^{-1}
\end{equation}
for any $\varepsilon\otimes 1\in\mathcal{O}_K^\times\otimes\F_p$.
By using the Chebotarev density theorem, we can choose 
primes $\frak{L}_1,\dots,\frak{L}_r$
of $K$ with properties (i), (ii), and
\begin{equation}\label{artin}
(\frak{L}_i,K(\sqrt[p]{\mathcal{O}_K^\times})/K)=\gamma_i\ \ (1\le i\le r).
\end{equation}
Now we choose rational integers $g_i$ such that
\begin{equation*}
g_i^{\frac{N(\frak{L_i})-1}{p}}\equiv \iota(1_{\F_p})\pmod{\frak{L}_i}.
\end{equation*}
Then $(g_i\mod{\frak{L}_i})\otimes 1$
generates $(\mathcal{O}_K/\frak{L}_i)^\times\otimes\F_p$,
and we obtain
\[
\phi_i(\sigma^{-1}(\varepsilon)\otimes 1)
=\gamma_i(\sqrt[p]{\sigma^{-1}(\varepsilon)})(\sqrt[p]
{\sigma^{-1}(\varepsilon)})^{-1}
\equiv
\sigma^{-1}(\varepsilon)^{\frac{N(\frak{L_i})-1}{p}}
\pmod{\frak{L}_i}
\]
by \eqref{gamma} and \eqref{artin}.
Hence if we assume that
\[
\mathrm{Pr}_i\circ f(\varepsilon\otimes 1)
=\sum_{\sigma\in \Delta\times G}a_\sigma^{(i)}\sigma,
\]
and 
\[(\sigma^{-1}(\varepsilon)\mod{\frak{L}_i})\otimes 1
=(g_i\mod{\frak{L}_i})\otimes b_\sigma^{(i)}
\]
with $a_\sigma^{(i)},\ b_\sigma^{(i)}\in\F_p$
for $1\le i\le r$ and $\sigma\in \Delta\times G$,
then
\begin{align*}
\iota(a_\sigma^{(i)})=
\phi_i(\sigma^{-1}(\varepsilon)\otimes 1)
&\equiv
\sigma^{-1}(\varepsilon)^{\frac{N(\frak{L_i})-1}{p}}
\\
&\equiv 
g_i^{b_\sigma\frac{N(\frak{L_i})-1}{p}}
\equiv\iota(b_\sigma^{(i)})
\pmod{\frak{L}_i},
\end{align*}
which implies
\[
a_\sigma^{(i)}=b_\sigma^{(i)},\ \ 
(\sigma^{-1}(\varepsilon)\mod{\frak{L}_i})\otimes 1
=(g_i\mod{\frak{L}_i})\otimes a_\sigma^{(i)}
\]
for all $i$ and $\sigma$.
Therefore we see that
\begin{equation*}
\begin{aligned}
\pi(\varepsilon\otimes 1)&=
(((\varepsilon\mod{\sigma\frak{L}_i})\otimes 1)_\sigma)_i
\\&=
(((g_i\mod{\sigma\frak{L}_i})\otimes a_\sigma^{(i)})_\sigma)_i
=\varphi\circ f(\varepsilon\otimes 1)
\end{aligned}
\end{equation*}
for any $\varepsilon\otimes 1\in\mathcal{O}_K^\times\otimes\F_p$.
This shows that the primes $\frak{L_1},\dots,\frak{L}_r$
and the rational integers $g_1,\dots,g_r$
also have property (iii).
\hfill$\Box$
%
%
%
%
%
\section{crucial proposition.}
Our aim in this section is to give the following, which plays
a crucial role in our construction:
\begin{prop}\label{prop1}
Let $p$ be a prime number,
and let $k$ be a totally imaginary number field of finite degree
with finite $\t{G}_k(p)$. 
Assume that
\[
r_2(k)\ge B_p(k).
\]
Then there exists a cyclic extension $k'/k$ of degree $p$ such that\\
(i)\ \ $k'\cap L_p(k)=k$, 
$L_p^{\mathrm{ab}}(Kk')=L_p^{\mathrm{ab}}(K)k'$
if we put $K:=L_p(k)(\mu_p)$,
hence $N_{Kk'/K}:\Cl(Kk')\otimes\Z_p\simeq\Cl(K)\otimes\Z_p$,\\
(ii)\ \ $L_p(k')=L_p(k)k'$, hence $\t{G}_{k'}(p)\simeq\t{G}_k(p)$
and $Kk'=L_p(k')(\mu_p)$,\\
(iii)\ \ $B_p(k')=B_p(k)$.

\end{prop}
We need some preliminary lemmas.
\begin{lem}\label{lem4}
Let $\frak{m}$ be an integral ideal of $k$.
Denote by $I_{k,\frak{m}}$ and $P_{k,\frak{m}}$ the group of the
 fractional ideals of $k$
which is prime to $\frak{m}$ and the subgroup of $I_{k,\frak{m}}$
which consists of all the principal ideals in $I_{k,\frak{m}}$, respectively.
Also we define $S_{k,\frak{m}}$ to be the subgroup of $P_{k,\frak{m}}$
which consists of all the principal ideals generated by elements congruent to 1
modulo $\frak{m}$.
Let $\frak{a}_1,\dots\frak{a}_{d(G)}$ be ideals of $\mathcal{O}_k$
which is prime to $\frak{m}$ and whose ideal classes
form a basis of the $p$-torsions $\Cl(k)[p]
:=\ker(\Cl(k)\overset{p}{\longrightarrow}\Cl(k))$ of the ideal class group
of $k$.
For $\alpha_i\in\mathcal{O}_k$ with $\frak{a}_i^p=\alpha_i\mathcal{O}_k$,
assume that there exist $\beta_i\in\mathcal{O}_k$ and
$\varepsilon_i\in\mathcal{O}_k^\times$ such that
\[
\alpha_i\equiv \beta_i^p\varepsilon_i\pmod{\frak{m}} 
\]
for $1\le i\le d(G)$.
Then we have
\[
\dim_{\F_p}((I_{k,\frak{m}}/S_{k,\frak{m}})\otimes\F_p)=
\dim_{\F_p}(\Cl(k)\otimes\F_p)+\dim_{\F_p}((P_{k,\frak{m}}/S_{k,\frak{m}})\otimes\F_p).
\]
\end{lem}
\noindent
{\bf Proof.}\ \ \ 
Since $\Cl(k)\simeq I_\frak{m}/P_\frak{m}$, we have the exact sequence
\[
0\longrightarrow
(P_m/S_m)[p]
\longrightarrow
(I_\frak{m}/S_\frak{m})[p]
\overset{\pi}{\longrightarrow}
\Cl(k)[p],
\]
where we denote by $A[p]$ the $p$-torsions of a module $A$.
To prove the lemma, it is enough to show the surjectivity of $\pi$.
Because $(\beta_i^{-1}\frak{a}_i)^p=\beta_i^{-p}\alpha_i\mathcal{O}_k$
and $\beta_i^{-p}\alpha_i\equiv \varepsilon_i\pmod{\frak{m}}
$ for each $i$,
we see that
\[
(\beta_i^{-1}\frak{a}_i\!\!\mod S_\frak{m})\in(I_\frak{m}/S_\frak{m})[p],\ 
\pi(\beta_i^{-1}\frak{a}_i\!\!\mod S_m)=(\frak{a}_i\!\!\mod P_m),
\]
which implies the surjectivity of $\pi$
because of the choice of the ideals $\frak{a}_i$.
\hfill$\Box$
\begin{lem}\label{lem5}
Put $K:=L_p(k)(\mu_p)$ and $G:=\Gal(L_p(k)/k)$.
Define $L/K$ to be the maximal elementary abelian $p$-subextension of
$L_p^{\ab}(K)/K$ such that $L/k(\m_p)$ is a Galois extension
and $G=\Gal(K/k(\mu_p))\simeq\Gal(K/k(\mu_p))$ acts on $\Gal(L/K)$
trivially via the inner automorphisms of $\Gal(L/k(\mu_p))$.
Also let $L'/K$ be the maximal elementary abelian $p$-subextension of
$L_p^{\ab}(K)/K$.
We assume that $\frak{L}_1,\dots,\frak{L}_r$ are primes of $K$ and
$k'/k$ a cyclic extension of degree $p$ which is linearly disjoint from
$L_p(k)/k$
such that\\
(i)\ \ 
$\Gal(L/K)=\langle
(\sigma\frak{L}_i,L/K)|\ 1\le i\le r,\ \sigma\in\Gal(K/k)
\rangle,
$\\
(ii)\ \ $L'k'$ is the maximal elementary abelian $p$-extension over $K$
which is unramified outside
$\{\sigma\frak{L}_i|1\le i\le r,\ \sigma\in\Gal(K/k)\}$,\\
(iii)\ \ if we denote by $\frak{l}_i$ the prime of $k$ below
$\frak{L}_i$,
then all of the primes $\frak{l}_1,\dots,\frak{l}_r$ ramify in $k'/k$.

\noindent
Then
$L_p^{\ab}(Kk')=L_p^{\ab}(K)k'$ and the norm map $N_{Kk'/K}$ induces
the isomorphism
\[
\Cl(Kk')\otimes\Z_p\overset{\sim}{\longrightarrow}\Cl(K)\otimes\Z_p.
\]
\end{lem}
\noindent
{\bf Proof.}\ \ \ 
It follows from assumption (i) and the isomorphism
\[
\Gal(L/K)\simeq\Gal(L^{\ab}_p(K)/K)_G\otimes\F_p
\]
 that
\[
\Gal(L_p^{\ab}(K)/K)=
\left\langle(\sigma\frak{L}_i,L^{\ab}(K)/K)
|1\le i\le r,\ \sigma\in\Gal(K/k)\right\rangle
\]
by Nakayama's lemma.
Hence if we denote by $\overline{\frak{L}}_i$ the prime
of $Kk'$ lying over $\frak{L}_i$, which is totally ramified in $Kk'/K$
by assumption (iii), then we have
\begin{equation}\label{ram}
\frak{L}_i=\overline{\frak{L}}_i^p
\end{equation}
and 
\begin{equation}\label{gen}
\Gal(L^{\ab}(K)k'/Kk')\!=\!
\left\langle(\sigma\overline{\frak{L}}_i,L^{\ab}(K)k'/Kk')
|1\!\le\! i\!\le\! r,\,\sigma\!\in\!\Gal(Kk'/k)\right\rangle\!.
\end{equation}

Let $M$ be the genus $p$-class field of the cyclic $p$-extension
$Kk'/K$, that is,
the maximal intermediate field of $L_p^{\ab}(Kk')/Kk'$
which is abelian over $K$.
We denote by $N/K$ the maximal abelian $p$-extension unramified
outside $\{\sigma\frak{L}_i|1\le i\le r,\ \sigma\in\Gal(K/k)\}$.
Then we have the inclusions
\[
K\subseteq Kk'\subseteq L_p^{\ab}(K)k'\subseteq M\subseteq N.
\]
It follows from assumption (ii) that
\begin{equation*}
\begin{aligned}
\dim_{\F_p}\Gal&(N/K)[p]=
\dim_{\F_p}\Gal(N/K)\otimes\F_p=\dim_{\F_p}\Gal(L'k'/K)\\
&=\dim_{\F_p}\Gal(L'/K)+1
=\dim_{\F_p}\Gal(L_p^{\ab}(K)/K)[p]+1,
\end{aligned}
\end{equation*}
which implies the cyclicity of $\Gal(N/L_p^{\ab}(K))$
by the exact sequence
\[
0\longrightarrow \Gal(N/L_p^{\ab}(K))[p]
\longrightarrow \Gal(N/K)[p]\longrightarrow
\Gal(L_p^{\ab}(K)/K)[p].
\]
Hence the primes of $L_p^{\ab}(K)$ lying over $\frak{L}_i$
totally ramify in $N/L_p^{\ab}(K)$
by assumption (iii) and the inclusions $L_p^{\ab}(K)\subseteq
L_p^{\ab}(K)k'\subseteq N$.
Therefore $N/L^{\ab}_p(K)k'$ has no non-trivial unramified subextensions,
hence $M=L^{\ab}_p(K)k'$.
Thus we find that
\begin{equation}\label{gcoinv}
\Gal(L^{\ab}_p(Kk')/Kk')_{\Gal(Kk'/K)}
\simeq
\Gal(M/Kk')=\Gal(L^{\ab}_p(K)k'/Kk').
\end{equation}
It follows from \eqref{gen}, the above isomorphism, and Nakayama's
lemma that
\begin{equation}\label{genn}
\Gal(L_p(Kk')/Kk')\!=\!
\langle
(\sigma\overline{\frak{L}}_i,L^{\ab}(Kk')/Kk')
|1\!\le\! i\!\le\! r,\,\sigma\in\Gal(Kk'/k)\rangle.
\end{equation}
Since $\sigma\overline{\frak{L}}_i\ (\sigma\in\Gal(Kk'/k))$
is a $\Gal(Kk'/K)$-invariant prime of $Kk'$ by \eqref{ram}, 
we see  that $(\sigma\overline{\frak{L}}_i, L_p^{\ab}(Kk')/Kk')$ is also 
$\Gal(Kk'/K)$-invariant
and derive from \eqref{gcoinv} and \eqref{genn} that
\begin{equation*}
\begin{aligned}
\Gal(L_p^{\ab}(K)k'/Kk')&\simeq
\Gal(L_p^{\ab}(Kk')/Kk')_{\Gal(Kk'/K)}\\
&=\Gal(L_p^{\ab}(Kk')/Kk'),
\end{aligned}
\end{equation*}
which implies
$
L_p^{\ab}(Kk')=L_p^{\ab}(K)k'.
$
Remaining assertion on the ideal class groups of $Kk'$ and $K$
follows from this identity and the functoriality of the Artin maps.
\hfill$\Box$

\

Now we shall give a proof of Proposition 1.\\
{\bf Proof of Proposition 1.}\ \ \ 
We first show that properties (ii) and (iii) of $k'/k$ follows from (i).
Assume that a cyclic extension $k'/k$ of degree $p$
satisfies (i).
If we put $K:=L_p(k)(\mu_p)$,
$G:=\Gal(L_p(k)/k)\simeq\Gal(K/k(\mu_p))$
and
\[
\Delta:=\Gal(k(\mu_p)/k)\simeq\Gal(K/L_p(k))\simeq\Gal(Kk'/L_p(k)k')
\]
as above,
we have 
\begin{equation*}
\begin{aligned}
\Gal(L_p^{\ab}(L_p(k)k')&/L_p(k)k')\simeq
\Gal(L_p^{\ab}(Kk')/Kk')_\Delta\\
&=\Gal(L_p^{\ab}(K)k'/Kk')_\Delta
\simeq
\Gal(L_p^{\ab}(K)/K)_\Delta\\
&\simeq\Gal(L_p^{\ab}(L_p(k))/L_p(k))=1,
\end{aligned}
\end{equation*}
hence $L_p^{\ab}(L_p(k)k')=L_p(k)k'$, which in turn implies
property (ii), namely, $L_p(k)k'=L_p(k')$ and $\t{G}_k(p)\simeq\t{G}_{k'}(p)$
because $L_p(k)k'/k'$ is unramified.
Then (i) also implies (iii) since 
$\t{G}_k(p)\simeq\t{G}_{k'}(p)$ and 
$
\Cl(L_p(k')(\mu_p))\otimes\Z_p=\Cl(Kk')\otimes\Z_p\simeq\Cl(K)\otimes\Z_p
$
as $G$-modules.
Thus it is enough to show that there exists $k'/k$ with property (i).

We write $L/K$ for the maximal elementary abelian $p$-subextension of
$L_p^{\ab}(K)/K$ such that $L/k(\m_p)$ is a Galois extension
and $G$ acts on $\Gal(L/K)$
trivially via the inner automorphisms of $\Gal(L/k(\mu_p))$
as in Lemma \ref{lem5}.
Then there exist $H\subseteq(\mathcal{O}_K^\times\otimes\F_p)^G$
and an abelian $p$-extension $M/K$ such that
\[
L=K(\sqrt[p]{H})M
\]
and $M$ is linearly disjoint from
$K(\sqrt[p]{\mathcal{O}_K^\times})/K$.
Here we note that
\begin{equation}\label{dimH}
\dim_{\F_p}\Gal(M/K),\ \dim_{\F_p}H\le
\dim_{\F_p}\Gal(L/K)=
\dim_{\F_p}(\Cl(K)\otimes\F_p)_G.
\end{equation}

Let $\frak{a}_1,\dots\frak{a}_{d(G)}$ be ideals of $\mathcal{O}_k$
whose ideal classes
form a basis of $\Cl(k)[p]$, and let $\alpha_i\in\mathcal{O}_k$ 
be integer such that $\frak{a}_i^p=\alpha_i\mathcal{O}_k$
as in Lemma \ref{lem4}.
It follows from the principal ideal theorem that
$\frak{a}_i\mathcal{O}_{L_p(k)}=\gamma_i\mathcal{O}_{L_p(k)}$ for some
$\gamma_i\in\mathcal{O}_{L_p(k)}$.
Hence there exists $\eta_i\in\mathcal{O}^\times_{L_p(k)}$ such that
\begin{equation}\label{eta}
\alpha_i=\eta_i\gamma_i^p
\end{equation}
for each $1\le i\le d(G)$.
We define $H_0\subseteq\mathcal{O}_{L_p(k)}^\times\otimes\F_p
=(\mathcal{O}_K^\times\otimes\F_p)^\bold{1}$
to be the $\F_p[G]$-submodule generated by $\eta_i\otimes 1$'s
$(1 \le i \le d(G))$.
We note that
\begin{equation}\label{imH_0}
H_0\subseteq((\mathcal{O}_K^\times\otimes\F_p)^\bold{1})^G
\end{equation}
%
%
%
by \eqref{eta}.
For each $\chi\in\widehat{\Delta}$,
let $F_\chi$ and $N_\chi$ be direct sum factors
of $(\mathcal{O}_K^\times\otimes\F_p)^\chi$ as given in Lemma \ref{lem2}.
We choose direct sum decomposition $F_\chi=F^0_\chi\oplus F^1_\chi$
of $F_\chi$ to free $\F_p[G]$-submodules
such that 
\begin{equation}\label{H^chi}
H^\chi\subseteq F^0_\chi\oplus N_\chi
\end{equation}
for each $\chi\in\widehat{\Delta}$,
%
%
where $H^\chi$ is the $\chi$-part of $H$.
It follows from facts 
$H\subseteq(\mathcal{O}_K^\times\otimes\F_p)^G$
and \eqref{dimH}
that we can choose $F^0_\chi$ with
\begin{equation}\label{dimF^0}
\rank_{\F_p[G]}F^0_\chi\le\dim_{\F_p}(\Cl(K)\otimes\F_p)_G.
\end{equation}
Under the above choice of $F^0_\chi$
we obtain the following estimation of $\rank_{\F_p[G]}F^1_\chi$
for any $\chi\in\widehat{\Delta}$ by using Lemma \ref{lem2}:
\begin{equation}\label{rkF^1pre}
\begin{aligned}
\rank_{\F_p[G]}F^1_\chi
\ge r_2(k)&-(\#G-1)(d(G)+r(G)+\dim_{\F_p}(\Cl(K)\otimes\F_p)_G)
\\&-\dim_{\F_p}(\Cl(K)\otimes\F_p)_G-1.
\end{aligned}
\end{equation}

Put 
\begin{equation*}
\begin{aligned}
m_0:&=([k(\mu_p):k](\#G-1)+1)(d(G)+r(G)+\dim_{\F_p}(\Cl(K)\otimes\F_p)_G)
\\&+\dim_{\F_p}(\Cl(K)\otimes\F_p)_G+1,\\
n_0:&=\dim_{\F_p}(\Cl(K)\otimes\F_p)_G.
\end{aligned}
\end{equation*}
Then it follows from Lemma \ref{lem2}, \eqref{dimH} and \eqref{dimF^0}
that
\begin{equation}\label{m_0}
m_0\ge\rank_{\F_p[G]}F_\chi^0+\dim_{\F_p}N_\chi^G+1,
\end{equation} 
\begin{equation}\label{n_0}
n_0\ge \dim_{\F_p}\Gal(M/K).
\end{equation}
Also we have
\begin{equation}\label{rkF^1}
\rank_{\F_p[G]}F_\chi^1\ge 2m_0+n_0+d(G)
\end{equation}
by the assumption $r_2(k)\ge B_p(k)$ and \eqref{rkF^1pre}.

\

For each $\chi\in\widehat{\Delta}$, we choose an $\F_p[G]$-homomorphism
\[
f^\chi:(\mathcal{O}_K^\times\otimes\F_p)^\chi
=F_\chi^0\oplus F_\chi^1\oplus N_\chi\longrightarrow
\bigoplus_{i=1}^{m_0+n_0}\F_p[G]\boldsymbol{a}_i^\chi,
\]
where $\boldsymbol{a}_i^\chi$'s are free $\F_p[G]$-basis, 
as follows:
\\
In the case where $\chi\ne \boldsymbol{1}$, we choose $f^\chi$ so that
\begin{equation}\label{chi-inj}
f^\chi|_{F_\chi^0\oplus N_\chi}:
F_\chi^0\oplus N_\chi\hookrightarrow 
\bigoplus_{i=1}^{m_0}\F_p[G]\boldsymbol{a}_i^\chi
\subseteq
\bigoplus_{i=1}^{m_0+n_0}\F_p[G]\boldsymbol{a}_i^\chi
\end{equation}
and
\begin{equation}\label{chi-surj}
f^\chi(F^1_\chi)=\bigoplus_{i=1}^{m_0+n_0}\F_p[G]\boldsymbol{a}_i^\chi.
\end{equation}
Such a homomorphism $f^\chi$ certainly exists by \eqref{m_0} and \eqref{rkF^1}.

\noindent
In the case where $\chi=\boldsymbol{1}$, we chose $f^{\boldsymbol{1}}$ so that
\begin{equation}\label{1-inj}
f^{\boldsymbol{1}}|_{F_{\boldsymbol{1}}^0\oplus N_{\boldsymbol{1}}}
:F_{\boldsymbol{1}}^0\oplus N_{\boldsymbol{1}}\hookrightarrow 
C:=\left\{\sum_{i=1}^{m_0}c_i\boldsymbol{a}_i^{\boldsymbol{1}}\left|\sum_{i=1}^{m_0}c_i=0\right.\right\}
\subseteq\bigoplus_{i=1}^{m_0+n_0}\F_p[G]\boldsymbol{a}_i^{\boldsymbol{1}}
\end{equation}
and, for some direct sum decomposition
\[
F_{\boldsymbol{1}}^1=F_{\boldsymbol{1}}^{10}
\oplus F_{\boldsymbol{1}}^{11}
\]
of $F_{\boldsymbol{1}}^1$
to free $\F_p[G]$-modules $F_{\boldsymbol{1}}^{10}$ and $F_{\boldsymbol{1}}^{11}$
with
$\rank_{\F_p[G]}F_{\boldsymbol{1}}^{10}=m_0-1$,
whose existence is assured by \eqref{rkF^1},
\begin{equation}\label{C-isom}
f^{\boldsymbol{1}}|_{F_{\boldsymbol{1}}^{10}}
:F_{\boldsymbol{1}}^{10}\simeq
C\subseteq\bigoplus_{i=1}^{m_0}\F_p[G]\boldsymbol{a}_i^{\boldsymbol{1}}
\end{equation}
and
\begin{equation}\label{D-surj}
f^{\boldsymbol{1}}(F_{\boldsymbol{1}}^{11})
=D:=
\left\{\sum_{i=1}^{m_0+n_0}c_i\boldsymbol{a}_i^{\boldsymbol{1}}
\left|\sum_{i=1}^{m_0+n_0}\varepsilon(c_i)=0
\right.\right\}
\subseteq\bigoplus_{i=1}^{m_0+n_0}\F_p[G]\boldsymbol{a}_i^{\boldsymbol{1}},
\end{equation}
where $\varepsilon:\F_p[G]\longrightarrow\F_p$ denotes the augmentation map.
Such a homomorphism $f^{\boldsymbol{1}}$ also certainly exists
by \eqref{m_0} and \eqref{rkF^1};
The fact $C\simeq\F_p[G]^{\oplus m_0-1}$
and \eqref{m_0} assure the existence of $f^{\boldsymbol{1}}$ with
properties \eqref{1-inj} and \eqref{C-isom}.
Also \eqref{rkF^1} and the fact
\[
D_G\simeq \F_p^{\oplus m_0+n_0+d(G)-1},
\]
which follows from the exact homology sequence
\begin{equation*}
\begin{aligned}
0=H_1(G,\F_p[G]^{\oplus m_0+n_0})
&\longrightarrow
H_1(G,\F_p)
\longrightarrow D_G\\
&\longrightarrow
\F_p^{\oplus m_0+n_0}
\longrightarrow \F_p\longrightarrow 0
\end{aligned}
\end{equation*}
and $H_1(G,\F_p)\simeq \F_p^{\oplus d(G)}$,
assure the existence of $f^{\boldsymbol{1}}$ with property \eqref{D-surj}.

It follows from Lemma \ref{lem3} and \eqref{n_0} that there exist
degree one prime ideals $\frak{L}_i$ of $K$ $(1\le i\le m_0+n_0)$, which lie
over distinct rational primes, and 
$g_i\in\Z$ such that\\
(a)\ \ \ $(\frak{L}_i,M/K)=1$ for $1\le i\le m_0$ and
$\langle(\frak{L}_i,M/K)|m_0+1\le i\le m_0+n_0\rangle=\Gal(M/K)$,\\
(b)\ \ \ $(g_i\mod\frak{L}_i)\otimes 1$ generates
$(\mathcal{O}_K/\frak{L}_i)^\times\otimes\F_p$, and
if we define the $\F_p[\Delta\times G]$
-isomorphism
\[
\varphi:\bigoplus_{i=1}^{m_0+n_0}\F_p[\Delta\times G]\t{\boldsymbol{a}}_i
\longrightarrow
\bigoplus_{i=1}^{m_0+n_0}\bigoplus_{\sigma\in \Delta\times G}
(\mathcal{O}_K/\sigma\frak{L}_i)^\times\otimes\F_p,
\]
where $\t{\boldsymbol{a}}_i=\sum_{\chi\in\widehat{\Delta}}
\boldsymbol{a}_i^\chi$ and $\Delta$ acts 
on $\boldsymbol{a}_i^\chi$ via $\chi$,
by
$\varphi(\t{\boldsymbol{a}}_i)=(g_i \mod \frak{L}_i)\otimes 1$,
and denote by
\[
\pi:\mathcal{O}_K^\times\otimes\F_p\longrightarrow
\bigoplus_{i=1}^{m_0+n_0}
\bigoplus_{\sigma\in \Delta\times G}
(\mathcal{O}_K/\sigma\frak{L}_i)^\times\otimes\F_p
\]
the natural projection map, then we have
\[
\pi=\varphi\circ(\oplus_{\chi\in\widehat{\Delta}}f^\chi),
\]
where
\[
\oplus_{\chi\in\widehat{\Delta}}f^\chi
:\mathcal{O}_K^\times\otimes\F_p\longrightarrow
\bigoplus_{\chi\in\widehat{\Delta}}
\bigoplus_{i=1}^{m_0+n_0}\F_p[G]\boldsymbol{a}_i^\chi
=\bigoplus_{i=1}^{m_0+n_0}\F_p[\Delta\times G]\t{\boldsymbol{a}}_i
\]
is the direct sum of $f^\chi$'s.

\

It follows from property (b) of $\frak{L}_i$'s,
\eqref{H^chi},
\eqref{chi-inj}, and \eqref{1-inj}
that 
\[
\pi|_H:H\hookrightarrow \bigoplus_{i=1}^{m_0}
\bigoplus_{\sigma\in \Delta\times G}
(\mathcal{O}_K/\sigma\frak{L}_i)^\times\otimes\F_p
\subseteq
\bigoplus_{i=1}^{m_0+n_0}
\bigoplus_{\sigma\in \Delta\times G}
(\mathcal{O}_K/\sigma\frak{L}_i)^\times\otimes\F_p,
\]
which implies that
\[
\Gal(K(\sqrt[p]{H})/K)
=\langle(\sigma\frak{L}_i, K(\sqrt[p]{H})/K)|\sigma\in\Delta\times G,\ 
1\le i\le m_0\rangle,
\]
and
$(\sigma\frak{L}_i, K(\sqrt[p]{H})/K)=1$ for $\sigma\in \Delta\times G$
and $m_0+1\le i\le m_0+n_0$.
Combining this and property (a) of $\frak{L}_i$'s,
we conclude that 
\begin{equation}\label{gen2}
\Gal(L/K)=\langle(\sigma\frak{L}_i, L/K)|\sigma\in \Delta\times G,\ 
1\le i\le m_0+n_0\rangle.
\end{equation}
We also derive from property (b),
\eqref{chi-surj}, \eqref{1-inj}, \eqref{C-isom} and \eqref{D-surj}
that
\[
\coker\,\pi=(\coker\,\pi)^{\Delta}\simeq\F_p,
\]
hence if we denote by $\frak{l}_i$
the prime of $k$ below $\frak{L}_i$ and 
by $P_{F,\frak{m}}/S_{F,\frak{m}}$
the ray class group of modulo
$\frak{m}:=\prod_{i=1}^{m_0+n_0}\frak{l}_i$
of a number field $F$ containing $k$
we have
\begin{equation}\label{dimstK}
(P_{K,\frak{m}}/S_{K,\frak{m}})\otimes\F_p
\simeq\F_p.
\end{equation}
We observe the maps
\begin{equation*}
\begin{aligned}
\mathcal{O}_k^\times\otimes\F_p\overset{\iota}{\longrightarrow}
((\mathcal{O}_K^\times\otimes\F_p)^{\boldsymbol{1}})^G
&\overset{\pi^{\boldsymbol{1},G}}{\longrightarrow}
\left(\bigoplus_{i=1}^{m_0+n_0}\bigoplus_{\sigma\in \Delta\times G}
(\mathcal{O}_K/\sigma\frak{L}_i)^\times\otimes\F_p\right)^{\Delta\times G}\\
&\simeq\bigoplus_{i=1}^{m_0+n_0}
(\mathcal{O}_k/\frak{l}_i)^\times\otimes\F_p\\
&\simeq(\mathcal{O}_k/\frak{m})^\times\otimes\F_p,
\end{aligned}
\end{equation*}
where $\iota$ is the natural map induced by the inclusion
$\mathcal{O}_k^\times\subseteq\mathcal{O}_K^\times$ and
$\pi^{\boldsymbol{1},G}$ is the restriction of $\pi$
to $((\mathcal{O}_K^\times\otimes\F_p)^{\boldsymbol{1}})^G$.

It follows from \eqref{1-inj}, \eqref{C-isom}, \eqref{D-surj}
and the fact $C\subseteq D$ that
\begin{equation}\label{im-pi'}
\begin{aligned}
&\im\,\pi^{\boldsymbol{1},G}
=
\varphi\circ
 f^{\boldsymbol{1}}(((\mathcal{O}_K^\times\otimes\F_p)^{\boldsymbol{1}})^G)
=\varphi\circ f^{\boldsymbol{1}}((F_{\boldsymbol{1}}^{10})^G)
+\varphi\circ f^{\boldsymbol{1}}((F_{\boldsymbol{1}}^{11})^G)
\\
&=
\varphi\circ f^{\boldsymbol{1}}(N_GF_{\boldsymbol{1}}^{10})
+\varphi\circ f^{\boldsymbol{1}}(N_GF_{\boldsymbol{1}}^{11})
=\varphi(N_GC)+\varphi(N_GD)=\varphi(N_GD)\\
&=
\left\{\!((g_1\!\!\!\!\mod\frak{l}_1)\!\otimes\! x_1,\!\cdots\!,
(g_{m_0+n_0}\!\!\!\!\mod\frak{l}_{m_0+n_0})\!\otimes\! x_{m_0+n_0})
\!\left|\sum_{i=1}^{m_0+n_0}\!\! x_i=0\right\}\right.
\end{aligned}
\end{equation}
and $\coker\,
\pi^{\boldsymbol{1},G}\simeq\F_p$.
Furthermore, because $F_{\boldsymbol{1}}^{11}\subseteq
(\mathcal{O}^\times_K\otimes\F_p)^{\boldsymbol{1}}
=\mathcal{O}^\times_{L_p(k)}\otimes\F_p$
and
$\im\,\pi^{\boldsymbol{1},G}=\varphi(N_GD)
=\pi^{\boldsymbol{1},G}(N_GF_{\boldsymbol{1}}^{11})$
by \eqref{im-pi'},
we see that $N_GF_{\boldsymbol{1}}^{11}\subseteq \im\,\iota$
and $\im\, (\pi^{\boldsymbol{1},G}\circ \iota)=\im\, \pi^{\boldsymbol{1},G}$.
Hence we have 
\begin{equation}\label{dimstk}
(P_{k,\frak{m}}/S_{k,\frak{m}})\otimes\F_p\simeq
\coker\,(\pi^{\boldsymbol{1},G}\circ \iota)=
\coker\,\pi^{\boldsymbol{1},G}\simeq \F_p.
\end{equation}

On the other hand, since $\pi(H_0)\subseteq\im\,\pi^{\boldsymbol{1},G}=
\im\, (\pi^{\boldsymbol{1},G}\circ \iota)$ by \eqref{imH_0},
we find from \eqref{eta} that
there exists $\varepsilon_i\in\mathcal{O}_k^\times$
such that
\begin{equation*}
(\alpha_i\!\!\mod\frak{m})\otimes 1=
(\eta_i\!\!\mod\frak{m})\otimes 1
=(\varepsilon_i\!\!\mod\frak{m})\otimes 1
\end{equation*}
for $1\le i\le d(G)$ in $(((\mathcal{O}_K/\frak{m})^\times\otimes\F_p)^{\boldsymbol{1}})^G
=(\mathcal{O}_k/\frak{m})^\times\otimes\F_p$.
Then it follows from \eqref{dimstk} and Lemma \ref{lem4} that
\[
\dim_{\F_p}(I_{k,\frak{m}}/S_{k,\frak{m}})\otimes\F_p=
\dim_{\F_p}(\Cl(k)\otimes\F_p)+1.
\]
Therefore there exists a cyclic extension $k'/k$
of degree $p$ and conductor dividing $\frak{m}$ which is linearly disjoint from
$L_p^{\ab}(k)/k$.
We derive from \eqref{dimstK} that
$L'k'/K$ is the maximal elementary abelian $p$-extension
which is unramified outside
$\{\sigma\frak{L}_i|1\le i\le m_0+n_0,\ \sigma\in\Gal(K/k)\}$.
We see that
every prime $\frak{l}_i$ ramifies in $k'/k$
since the natural map
$(\mathcal{O}_k/\frak{l}_i)^\times\otimes\F_p\longrightarrow
\coker\,(\pi^{\boldsymbol{1},G}\circ \iota)$
is injective by \eqref{im-pi'},
whose image is isomorphic to the inertia subgroup for the prime $\frak{l}_i$
of the Galois group of the maximal elementary abelian $p$-extension
over $k$ unramified outside
$\{\frak{l}_1,\dots,\frak{l}_{m_0+n_0}\}$.
Hence, combining \eqref{gen2}, we find that
the extension $k'/k$ satisfies conditions (i)--(iii) of Lemma \ref{lem5}.
Thus we conclude by using Lemma \ref{lem5}
that $k'/k$ has property (i) of Proposition \ref{prop1},
which completes the proof.
\hfill $\Box$
\section{Embedding problem of Galois extensions}
In this section, we recall some facts from 
embedding problem of Galois extension
to prove Proposition II
(consult, for example, \cite{Neu} and \cite{Led} for this section).
\par
Let $p$ be any prime number and $F/E$ a finite $p$-extension
of number fields of finite degree
such that $G:=\Gal(F/E)$ fits the following exact sequence with a
finite $p$-group $G'$:
\begin{equation}\label{ep}
1 \longrightarrow \Z/p \overset{\iota}{\longrightarrow}G'
\overset{\pi}{\longrightarrow} G\longrightarrow 1.
\end{equation}
The Galois extension $M/E$ containing $F$
is called a solution of embedding problem \eqref{ep}
if and only if there exists an isomorphism $\theta:G'\simeq\Gal(M/E)$
such that the diagram
\begin{equation*}
\begin{CD}
G' @>{\pi}>> G\\
@V{\theta\ \wr}VV @|\\
\Gal(M/E) @>{\mathrm{restriction}}>> G
\end{CD}
\end{equation*}
is commutative.

We recall that every cyclic extension $M/F$ of degree $p$ is
given by the unique degree $p$ subextension of the cyclic extension
$F(\mu_p,\,\sqrt[p]{\alpha})/F$
for some non-trivial $\alpha\otimes 1\in(F(\mu_p)^\times\otimes\F_p)^\omega$,
$\omega$ being the cyclotomic character
$\Gal(F(\mu_p)/F)\longrightarrow\F_p^\times$ as above.
 In what follows, we denote by $F\{\sqrt[p]{\alpha}\}/F$
the unique subextension of $F(\mu_p,\,\sqrt[p]{\alpha})/F$
of degree $p$ for any non-trivial
$\alpha\otimes 1\in(F(\mu_p)^\times\otimes\F_p)^\omega$.
Now we recall the following theorem, which plays a crucial roll
to establish the final step of our construction:
\begin{thma}
Let $p$ be a prime number and $F/E$ an unramified
$p$-extension of number fields of
finite degree with Galois group $G$.
Then embedding problem \eqref{ep} always has a solution.
Furthermore, if $F\{\sqrt[p]{\alpha_0}\}/E$ is a solution for
$\alpha_0\otimes 1\in(F(\mu_p)^\times\otimes\F_p)^\omega$,
then all the solutions are given by
$F\{\sqrt[p]{a\alpha_0}\}/E$ for
$a\otimes 1\in(E(\mu_p)^\times\otimes\F_p)^\omega$
with $a\alpha_0\otimes 1\ne 0$ in $(F(\mu_p)^\times\otimes\F_p)^\omega$.
\end{thma}
\noindent{\bf Proof.}\ \ The existence of a solution follows from
\cite[Satz 2.2, 4.7, and 5.1]{Neu}, and the description of
all the solutions is given by \cite[Cor.8.1.5]{Led}.
\hfill $\Box$
\section{Final step of the construction}
We shall finish the proof of Proposition II in this section.
Our first aim is to show the following:
\begin{lem}\label{lemma6}
Let $p$ be a prime number
and $k$ a totally imaginary number field of finite degree.
Assume that $G:=\t{G}_k(p)$ is finite and
$r_2(k)\ge B_p(k)$.
Then for any given exact sequence of $p$-groups
\begin{equation}\label{gpext}
1\longrightarrow \Z/p
\longrightarrow G'\overset{\pi}\longrightarrow G
\longrightarrow 
1,
\end{equation}
there exists a finite extension $k'/k$ such that
\\
(i)\ \ $L_p(k')=L_p(k)k'$ and $k'\cap L_p(k)=k$,
hence the restriction induces $\t{G}_{k'}(p)\simeq\t{G}_k(p)$,\\
(ii)\ \ $B_p(k')=B_p(k)$,\\
(iii)\ \ there exists an 
$\varepsilon_0\otimes 1\in(\mathcal{O}_{L_p(k')(\mu_p)}^\times\otimes\F_p)^\omega$
such that 
\[
\begin{CD}
 G' @>{\ \ \ \ \pi\ \ \ \ }>> G\\
@V{\theta\ \wr}VV @AA{\wr\ \mathrm{restriction}}A\\
\Gal(L_p(k')\{\sqrt[p]{\varepsilon_0}\}/k')
@>{\mathrm{restriction}}>> \t{G}_{k'}(p)
\end{CD}
\]
is commutative with some isomorphism
$\theta:G'\simeq\Gal(L_p(k')\{\sqrt[p]{\varepsilon_0}\}/k')$.
\end{lem}
\noindent{\bf Proof.}\ \ 
First we assume that group extension \eqref{gpext} splits.
One can show that
\[
 \dim_{\F_p}(\mathcal{O}_{k(\mu_p)}^\times\otimes\F_p)^\omega
=r_2(k),\ \mbox{or}\ r_2(k)+1
\]
by the same manner as in the proof of Lemma \ref{chi-unit}. 
Then it follows from the assumption 
$r_2(k)\ge B_p(k)$
that 
\[
 \dim_{\F_p}(\mathcal{O}_{k(\mu_p)}^\times\otimes\F_p)^\omega
> d(G).
\]
Hence there exists 
$\varepsilon_0\otimes 1\in(\mathcal{O}_{k(\mu_p)}^\times\otimes\F_p)^\omega$
such that $k\{\sqrt[p]{\varepsilon_0}\}/k$ is a cyclic extension of
degree $p$ which is linearly disjoint from $L_p(k)/k$.
Then we have $\Gal(L_p(k)\{\sqrt[p]{\varepsilon_0}\}/k)\simeq G'$.
Thus $k'=k$ is a required extension.
\par
Now we assume that group extension \eqref{gpext} does not split.
It follows from Theorem A that
there exists a solution $M/k$ of embedding problem \eqref{gpext}.
Then there is an element $\alpha\otimes 1\in (K^\times\otimes\F_p)^\omega$,
$K$ denoting $L_p(k)(\mu_p)$ as above,
such that $M=L_p(k)\{\sqrt[p]{\alpha}\}$.
Since $M/k$ is a Galois extension, we see
$\alpha\otimes 1\in ((K^\times\otimes\F_p)^\omega)^G$.
Hence it follows from the exact cohomology sequence
\[
0\longrightarrow I_K^G\overset{p}{\longrightarrow}I_K^G
\longrightarrow (I_K\otimes\F_p)^G
\longrightarrow H^1(G,I_K)=0,
\]
$I_K$ being the group of the fractional ideals of $K$, 
that
\begin{equation}\label{Aa}
\alpha\mathcal{O}_{K}=\frak{A}^p\frak{a}
\end{equation}
for some ideals $\frak{A}$ of $K$ and $\frak{a}$ of $k(\mu_p)$.
We denote by $h$ the non-$p$-part of the class number of $K$,
which is divisible by that of $k(\mu_p)$.
Then by replacing $\alpha$ to $\alpha^h$,
we may assume that the orders of the ideal class
$\cl_K(\frak{A})\in\Cl(K)$
and $\cl_{k(\mu_p)}(\frak{a})\in\Cl(k(\mu_p))$
containing $\frak{A}$ and $\frak{a}$,
respectively, are powers of $p$.
\par
By using Proposition \ref{prop1} repeatedly,
we obtain a tower of extensions $k_n/k\ (n\ge 0)$ of degree $p^n$ with
$k_l\subseteq k_m$ ($l\le m$)
such that $L_p(k_n)=L_p(k)k_n$, $k_n\cap L_p(k)=k$, $\t{G}_{k_n}(p)\simeq
G$, and
\begin{equation}\label{nisom}
N_{m,n}:=N_{Kk_m/Kk_n}\colon\Cl(Kk_m)\otimes\Z_p\simeq\Cl(Kk_n)\otimes\Z_p
\end{equation}
for $m\ge n\ge 0$.
Denote by
$j_{n,m}:\Cl(k_n(\mu_p))\otimes\Z_p
\longrightarrow\Cl(k_m(\mu_p))\otimes\Z_p$ and
$J_{0,n}\colon\Cl(K)\otimes\Z_p\longrightarrow\Cl(Kk_n)\otimes\Z_p$
the natural homomorphisms induces by the inclusions
$k_n(\mu_p)\subseteq k_m(\mu_p)$ and $K\subseteq Kk_n$, respectively. 
It follows from \eqref{nisom} and $N_{n,0}\circ J_{0,n}=p^n$
that
\begin{equation}\label{kerJ}
\ker J_{0,n}=\Cl(K)[p^n].
\end{equation}
Also, since we obtain the estimation
\begin{equation*}
\begin{aligned}
\#(\Cl(k_n(\mu_p))\otimes\Z_p)
&\le [Kk_n:k_n(\mu_p)]\#(\Cl(Kk_n)\otimes\Z_p)_G
\\&=[K:k(\mu_p)]\#(\Cl(K)\otimes\Z_p)_G
\end{aligned}
\end{equation*}
for all $n\ge 0$ by using \eqref{nisom},
we find from the surjectivity of
the maps 
\[
N_{m,n}\colon
 \Cl(k_m(\mu_p))\otimes\Z_p\longrightarrow\Cl(k_{n}(\mu_p))\otimes\Z_p
\]
for $m\ge n\ge 0$
that there is a  number $n_0\ge 0$ such that
\[
N_{n,n_0}\colon
\Cl(k_n(\mu_p))\otimes\Z_p\longrightarrow\Cl(k_{n_0}(\mu_p))\otimes\Z_p
\]
is an isomorphism for any $n\ge n_0$.
Hence we see that 
\begin{equation}\label{kerj}
\Cl(k(\mu_p))[p^{n-n_0}]\subseteq\ker j_{0,n}
\end{equation}
for $n\ge n_0$ by $N_{n,n_0}\circ j_{n_0,n}=p^{n-n_0}$.
Therefore, by \eqref{kerJ} and \eqref{kerj},
we conclude that $\frak{A}$ and $\frak{a}$
are principal ideals in $Kk_n$ and $k_n(\mu_p)$, respectively,
if $n$ is sufficiently large:
\begin{equation}\label{principal}
\frak{A}\mathcal{O}_{Kk_n}=A\mathcal{O}_{Kk_n},\ \ 
\frak{a}\mathcal{O}_{k_n(\mu_p)}=a\mathcal{O}_{k_n(\mu_p)},\ \ 
A\in (Kk_n)^\times,\ \ a\in k_n(\mu_p)^\times.
\end{equation}

It follows from \eqref{Aa} and \eqref{principal} that there exists
$\varepsilon_0\otimes 1\in(\mathcal{O}_{Kk_n}^\times\otimes\F_p)^\omega$
such that
\[
 Mk_n=L_p(k_n)\{\sqrt[p]{\alpha}\}=
L_p(k_n)\{\sqrt[p]{\varepsilon_0 a^\omega}\},
\]
where $a^\omega\in k_n(\mu_p)^\times$
is so that
$a^\omega\otimes 1$ is the projection of
$a\otimes 1\in k_n(\mu_p)^\times\otimes\F_p$ to
$(k_n(\mu_p)^\times\otimes\F_p)^\omega$.
Here, $L_p(k_n)\{\sqrt[p]{\varepsilon_0 a^\omega}\}/k_n$
is a solution of embedding problem
\[
1\longrightarrow \Z/p\longrightarrow
G'\overset{r^{-1}\circ\pi}\longrightarrow 
\Gal(L_p(k_n)/k_n)\longrightarrow 1,
\]
where $r:\Gal(L_p(k_n)/k_n)\simeq G$ is the isomorphism
induced by the restriction.
Since the above group extension is not split from our assumption,
we see that $\varepsilon_0\otimes 1\ne 0$, which combined with Theorem A implies
$L_p(k_n)\{\sqrt[p]{\varepsilon_0}\}/k_n$ is also a solution of
this embedding problem.
Therefore $k':=k_n$ is a required extension of $k$.
\hfill$\Box$

\

We will employ the following lemma in
our final step of the construction
to assure a constructed unramified
$p$-extension to be maximal:
\begin{lem}\label{cent}
Let $p$ be any prime number, $F$ a number field of finite degree
with $L_p(F)=F$ and $S$ a finite set of primes of $F$. 
We denote by $F_S/F$ the maximal elementary abelian $p$-extension
unramified outside $S$. For any prime $v$ of $F$,
we write $D_v$ for the decomposition subgroup
of $\Gal(F_S/F)$ at $v$.
We assume that the map 
\[
\phi:\underset{v\,:\,\mbox{\tiny prime of $F$}}{\bigoplus}H_2(D_v,\Z)
\longrightarrow H_2(\Gal(F_S/F),\Z)
\]
induced by the natural inclusions $D_v\subseteq\Gal(F_S/F)$
is surjective. Then we have $L_p(F_S)=F_S$.
\end{lem}
%
%
\noindent{\bf Proof.}\ \ 
We first note that $F_S/F$
is a finite extension. Let $M$ be the genus $p$-class field of $F_S/F$,
namely, the maximal unramified $p$-extension $M$ of $F_S$ such that
there exists an abelian extension $F'/F$ with $M=F_SF'$.
Also, we denote by $L$ the central $p$-class field of $F_S/F$,
namely, $L/F$ is the maximal normal sub-extension
of $L_p(F_S)/F$ containing $F_S$ such that
$\Gal(L/F_S)$ is contained in the center of $\Gal(L/F)$.
In other words, $L$ is the fixed field of $L_p(F_S)$
by the commutator subgroup $(\Gal(L_p(F_S)/F_S),\Gal(L_p(F_S)/F))$
of $\Gal(L_p(F_S)/F_S)$.
Then $M$ is contained in $L$ and we have a surjection
\begin{equation*}
\mathrm{coker}\,\phi\longrightarrow\Gal(L/M)
\end{equation*}
(see \cite[Theorem 3.11]{Fro83} and section 11.4 in \cite{TA}).
 Hence we have $M=L$ by the assumption $\coker\,\phi=0$.
Suppose that $F_SF'/F_S$ is a non-trivial unramified $p$-extension
for an abelian extension $F'/F$.
Then it follows from the definition of $F_S$ that
$F'/F$ is unramified outside $S$, and that $F'$ must contain a cyclic extension
$F_0/F$ of degree $p^2$ over $F$, which is totally ramified
at some ramified prime by $L_p(F)=F$.
Hence $F_SF'/F_S$ is not unramified, which is a contradiction.
Thus we conclude that $L=M=F_S$.
\par
Since $\Gal(L_p(F_S)/F)$ is a pro-$p$-group, 
it follows from the isomorphism
\[
1=\Gal(L/F_S)\simeq\Gal(L_p(F_S)/F_S)
/(\Gal(L_p(F_S)/F_S),\Gal(L_p(F_S)/F))
\]
that
$L_p(F_S)=F_S$.
Therefore we obtain the lemma.
\hfill$\Box$
%
%
%
%
%

\

Now, by virtue of Lemma \ref{lemma6},
we may assume that $k$ is an totally imaginary number field
of finite degree such that $G:=\t{G}_k(p)$ is finite,
every prime of $k$ lying over $p$ splits completely in $L_p(k)$,
$r_2(k)\ge B_p(k)$, and there exists $\varepsilon_0\in
(\mathcal{O}_{K}^\times\otimes\F_p)^\omega\ \ (K:=L_p(k)(\mu_p))$
such that $L_p(k)\{\sqrt[p]{\varepsilon_0}\}/k$
is a solution of the embedding problem
\begin{equation}\label{eprob}
1\longrightarrow \Z/p\longrightarrow G'
\overset{\pi}{\longrightarrow}G\longrightarrow 1
\end{equation}
in the statement of Proposition II.
Furthermore, we may assume also that
\begin{equation}\label{B_p(k)}
\begin{aligned}
r_2(k)&\ge
\left((2(p-1)+1)\left(\#G'-1\right)
+2\right)\times\\
&\left(
d(G')\!+\!r(G')
\!+\!p^2(\dim_{\F_p}\!\left((\Cl\left(L_p(k)(\mu_p))\otimes\F_p\right)
+2(p-1)\#G)\right)\right)\\
&+d(G')
+4p^2(\dim_{\F_p}\!\left((\Cl\left(L_p(k)(\mu_p))\otimes\F_p\right)
+2(p-1)\#G\right)+3
\end{aligned}
\end{equation}
by replacing $k$ with some suitable finite extension of $k$
by using Proposition 1 repeatedly.
In what follows, we shall observe the global unit
$\varepsilon_0$ locally at $p$.

We introduce some notations which will be used in what follows;
For any number field $F$ of finite degree
and prime $\frak{p}$ of $F$ lying above $p$,
we denote by $U_\frak{p}(F)$ the pro-$p$-part of the local unit
group of the completed field $F_\frak{p}$ of $F$ at $\frak{p}$,
and define the pro-$p$-part of the semi-local unit group at $p$ by
$U(F)=\bigoplus_{\frak{p}|p}U_\frak{p}(F)$.
We embed the unit group of the localization $\mathcal{O}_{F,p}$ of 
the maximal order $\mathcal{O}_{F}$ of $F$ at $p$
diagonally into $U(F)$ as usual. 
In the case where $\mu_p\subseteq F$,
we define $U'_\frak{p}(F)$ to be the submodule of $U_\frak{p}(F)$ consisting
of all the elements $u\in U_\frak{p}(F)$ such that 
$F_\frak{p}(\sqrt[p]{u})/F_\frak{p}$ is an unramified extension
(including the case $F_\frak{p}(\sqrt[p]{u})=F_\frak{p}$) for
a prime $\frak{p}$ lying over $p$,
and put $U'(F)=\bigoplus_{\frak{p}|p}U'_\frak{p}(F)$.
Then we have the inclusions $U(F)^p\subseteq U'(F)\subseteq U(F)$
and put $R(F)=U(F)/U(F)^p$ and $R'(F)=U(F)/U'(F)$.
Here we note that
\begin{equation}\label{unru}
U'_\frak{p}(F)/U_\frak{p}(F)^p\simeq\F_p
\end{equation}
as $\F_p$-modules.

Since every primes of $k$ lying over $p$ splits completely in $L_p(k)/k$
by the assumption, 
we see that 
\begin{equation}\label{tensor}
R(K)\simeq R(k(\mu_p))\underset{\F_p}{\otimes}\F_p[G],\ 
R'(K)\simeq R'(k(\mu_p))\underset{\F_p}{\otimes}\F_p[G].
\end{equation}
For a prime of $\frak{p}$ of $k$ lying over $p$,
the $p$-adic logarithm map $\log_p$ and
the normal basis theorem induces
$\Gal(k(\mu_p)_{\overline{\frak{p}}}/k_{\frak{p}})$-isomorphisms
\[
U_{\overline{\frak{p}}}(k(\mu_p))\underset{\Z_p}{\otimes}\Q_p
\!\overset{\log_p\otimes\Q_p}{\simeq}\!
\mathcal{O}_{k(\mu_p)_{\overline{\frak{p}}}}\underset{\Z_p}{\otimes}\Q_p
\!\simeq\!
k(\mu_p)_{\overline{\frak{p}}}
\!\simeq\!
\Q_p[\Gal(k(\mu_p)_{\overline{\frak{p}}}/k_{\frak{p}})]
^{\oplus [k_{\frak{p}}:\Q_p]},
\]
where $\overline{\frak{p}}$ is a fixed prime of $k(\mu_p)$ lying over $\frak{p}$
and $\mathcal{O}_{k(\mu_p)_{\overline{\frak{p}}}}$
stands for the valuation ring of
$k(\mu_p)_{\overline{\frak{p}}}$.
Hence we deduce 
\[
R(k(\mu_p))^\omega\simeq\F_p^{\oplus[k:\Q]+s},\ \ 
R'(k(\mu_p))^\omega\simeq\F_p^{\oplus[k:\Q]}
\]
as $\F_p$-modules, where $s$ denotes the number of the primes of $k$
lying over $p$,
from the facts
\[
R(k(\mu_p))\simeq
\bigoplus_{\frak{p}|p}\left(
U_{\overline{\frak{p}}}(k(\mu_p))
\underset{\Z_p[\Gal(k(\mu_p)_{\overline{\frak{p}}}/k_{\frak{p}})]}
{\otimes}\Z_p[\Gal(k(\mu_p)/k)]
\right)\otimes\F_p,
\]
and
$U'_{\overline{\frak{p}}}(k(\mu_p))/U_{\overline{\frak{p}}}(k(\mu_p))^p
\subseteq(U_{\overline{\frak{p}}}(k(\mu_p))\otimes\F_p)
^{\omega|_{\Gal(k(\mu_p)_{\overline{\frak{p}}}/k_{\frak{p}})}}$,
and \eqref{unru}
by the similar manner as in the proof of Lemma \ref{chi-unit}.
Hence we have
\begin{equation}\label{rankR}
R(K)^\omega\simeq \F_p[G]^{\oplus[k:\Q]+s},\ \ 
R'(K)^\omega\simeq \F_p[G]^{\oplus[k:\Q]}
\end{equation}
as $G$-modules
by using \eqref{tensor}.

Since $L_p(k)\{\sqrt[p]{\varepsilon}\}/L_p(k)$ is an abelian $p$-extension
unramified outside $p$ 
for any $\varepsilon\otimes 1\in(\mathcal{O}_K^\times\otimes\F_p)^\omega$
and $L_p(k)$ has no non-trivial unramified $p$-extension,
the natural maps
$(\mathcal{O}_K^\times\otimes\F_p)^\omega\rightarrow R'(K)^\omega$
and $(\mathcal{O}_K^\times\otimes\F_p)^\omega\rightarrow R(K)^\omega$
is injective. We put $E:=\im((\mathcal{O}_K^\times\otimes\F_p)^\omega\rightarrow R(K)^\omega)$.
\begin{lem}\label{lem-8}
Let $Z$ be the kernel of the natural projection map
$\pi:R(K)^\omega\longrightarrow R'(K)^\omega$.
Then there exist free $\F_p[G]$-submodules $P$ and $Q$ of $R(K)^\omega$
such that 
\begin{equation*}
R(K)^\omega=P\oplus Q\oplus Z,\ 
R'(K)^\omega\overset{(\pi|_{P\oplus Q})^{-1}}
{\simeq}P\oplus Q,
\end{equation*}
and that $E\subseteq P$ and
$
\rank_{\F_p[G]}Q\geq
r_2(k)-1-([k(\mu_p):k](\#G-1)+1)(d(G)+r(G)+\dim_{\F_p}(\Cl(K)\otimes\F_p)_G).
$
\end{lem}
\noindent{\bf Proof.}\ \ 
Let $r=\dim_{\F_p}E^G$. Then we see that
there exists an injection of $\F_p[G]$-modules
$h:E\longrightarrow \F_p[G]^{\oplus r}$.
Let $i:E\rightarrow R(K)^\omega$
be the natural inclusion map.
Since $R(K)^\omega\simeq\F_p[G]^{\oplus [k:\Q]+s}$ by \eqref{rankR},
$R(K)^\omega$ is an injective object in the category of
finitely generated $\F_p[G]$-modules.
Hence there exists a homomorphism 
$g:\F_p[G]^{\oplus r}\longrightarrow R(K)^\omega$
such that $g\circ h=i$.
\par
We find that $g$ is injective as follows.
Let $y\in(\mathrm{ker}\,g)^G$ be any element.
 Since $E^G\simeq\F_p^{\oplus r}\simeq(\F_p[G]^{\oplus r})^G$,
the injection $h$ induces the isomorphism $E^G\simeq(\F_p[G]^{\oplus r})^G$.
Hence $y=h(x)$ for some $x\in E^G$ and we have $0=g(y)=i(x)$,
which implies $x=0$ and $y=0$. Thus we show that $(\mathrm{ker}\,g)^G=0$,
which in turn is equivalent to $\mathrm{ker}\,g=0$.
\par
Put $P=\mathrm{im}(g)$. Then $P\simeq\F_p[G]^{\oplus r}$ and $E\subseteq P$.
Let $x\in (P\cap Z)^G$ be any element. 
Because, as we have seen, $P^G=E^G$ and $E\cap Z=0$ holds,
we have $x=0$. Hence we deduce $P\cap Z=0$, then $P+Z$
is a direct sum in $R(K)^\omega$.

Because any direct factor of a free $\F_p[G]$-module is also free,
we deduce from \eqref{rankR} that 
\begin{equation}\label{rankN}
Z\simeq\F_p[G]^{\oplus s}.
\end{equation}
Hence $P\oplus Z$ is also a free $\F_p[G]$-module.
Then injectivity of $P\oplus Z$ shows that
there exists a free $\F_p[G]$-submodule
$Q$ of $R(K)^\omega$ such that $R(K)^\omega=P\oplus Q\oplus Z$.

Now we shall estimate the $\F_p[G]$-rank of $Q$.
It follows from Lemmas \ref{chi-unit} and \ref{lem2} that
\[
E\simeq F_\omega\oplus N_\omega
\]
for a certain free $\F_p[G]$-module $F_\omega$ of rank $r_\omega\le r_2(k)+1$
and $\F_p[G]$-module $N_\omega$ with
$\dim_{\F_p}N_\omega^G\le
([k(\mu_p):k](\# G-1)+1)(d(G)+r(G)+\dim_{\F_p}(\Cl(K)\otimes\F_p)_G)$.
Hence we have
\begin{equation}\label{rankP}
\begin{aligned}
&\rank_{\F_p[G]}P=\dim_{\F_p}E^G\le r_2(k)+1\\
&+([k(\mu_p):k](\# G-1)+1)(d(G)+r(G)+\dim_{\F_p}(\Cl(K)\otimes\F_p)_G).
\end{aligned}
\end{equation}
Therefore we derive from \eqref{rankR}, \eqref{rankN} and \eqref{rankP} that
\begin{equation*}
\begin{aligned}
&\rank_{\F_p[G]}Q=2r_2(k)-\rank_{\F_p[G]}P\ge r_2(k)-1\\
&-([k(\mu_p):k](\# G-1)+1)(d(G)+r(G)+\dim_{\F_p}(\Cl(K)\otimes\F_p)_G).
\end{aligned}
\end{equation*}
\hfill$\Box$
\begin{lem}\label{lem-9}
Let $F/K$ be any finite abelian $p$-extension 
linearly disjoint from the maximal abelian extension of $K$
which is unramified outside $p$.
Then for any $u\in R(K)$ and $\tau\in\Gal(F/K)$,
there exist infinitely many 
degree one principal prime ideals $\Lambda\mathcal{O}_{K}$
of $\mathcal{O}_{K}$ such that $\Lambda\mathcal{O}_{K}$
is prime to $p$, $(\Lambda\!\!\!\mod U(K)^p)=u$ in $R(K)$
and $(\Lambda\mathcal{O}_K, F/K)=\tau$.
\end{lem}
\noindent{\bf Proof.}\ \ 
Let $M$ be the maximal abelian $p$-extension of $K$
which is unramified outside $p$.
We obtain the exact sequence
\begin{equation}\label{eseq-K}
\mathcal{O}_K^\times\otimes\F_p\longrightarrow
R(K)
\overset{\rho}{\longrightarrow}
\Gal(M/L^{\ab}_p(K))\otimes\F_p\longrightarrow 0
\end{equation}
by class field theory,
where $\rho$ is the map induced by the reciprocity map
from the idele group of $K$ to $\Gal(M/K)$.
Let $H$ be the maximal unramified abelian extension of 
$K$ and $\sigma\otimes 1=\rho(u)\in\Gal(M/L_p^{\ab}(K))\otimes\F_p$.
Then, it follows from the Chebotarev density theorem
that there exist infinitely many degree one 
primes $\frak{L}$ of $K$ not lying over $p$
such that $(\frak{L},H/K)=1$, 
$(\frak{L},M/K)=\sigma^{-1}$ and 
$(\frak{L},F/K)=\tau$,
since $\sigma\in\Gal(M/L_p^{\ab}(K))$, $H\cap M=L_p^{\ab}(K)$,
and the fact that $F$ is linearly disjoint to $M$ and $H$
over $K$. 
The assumption $(\frak{L},H/K)=1$ implies that
$\frak{L}$ is a principal prime ideal, say, 
$\frak{L}=\Lambda_0\mathcal{O}_{K}$,
and we deduce from the assumption $(\frak{L},M/K)=\sigma^{-1}$
that $\rho(\Lambda_0\!\!\!\mod U(K)^p)=\sigma\otimes 1$.
Therefore it follows from \eqref{eseq-K} that
$(\Lambda \!\!\!\mod U(K)^p)=u$ in $R(K)$
and $(\Lambda\mathcal{O}_{K}, F/K)=(\frak{L},F/K)
=\tau$ for a prime element $\Lambda=\Lambda_0\varepsilon$ with some 
$\varepsilon\in\mathcal{O}_{K}^\times$.
\hfill$\Box$

\

Recall that $L_p(k)\{\sqrt[p]{\varepsilon_0}\}/k$ is a solution of
embedding problem \eqref{eprob}
with $\varepsilon_0\otimes 1\in(\mathcal{O}_K^\times\otimes\F_p)^\omega$.
Let $R(K)^\omega=P\oplus Q\oplus Z$ be a decomposition given
by Lemma \ref{lem-8}.
Since $L_p(k)\{\sqrt[p]{\varepsilon_0}\}/k$ is a Galois extension
and $P$ is a free $\F_p[G]$-module,
we have $(\varepsilon_0\mod U(K)^p)\in E^G\subseteq P^G=N_GP$.
Hence 
\begin{equation}\label{c1}
(\varepsilon_0\mod U(K)^p)=N_Gx
\end{equation}
for some $x\in P$.
\par
Let $\{\sigma_1,\sigma_2,\dots,\sigma_d\}$
be a generator system of $G$, where $d=d(G)$,
and $\{q_1,q_2,\dots,q_t\}$ ($t=\rank_{\F[G]}Q$) be a system of
free basis of the free $\F_p[G]$-module $Q$:
$Q=\bigoplus_{i=1}^t\F_p[G]q_i$.
Here we have 
\begin{equation}\label{td}
t\geq 2d
\end{equation}
by Lemma \ref{lem-8} and the assumption $r_2(k)\ge B_p(k)$.

By virtue of Lemma \ref{lem-9} and \eqref{td}, there exists
a degree one prime element
$\Lambda_1\in\mathcal{O}_{K}$ which is prime to $p$
such that
\begin{equation}\label{c2}
(\Lambda_1\mod U(K)^p)=-x+\sum_{i=1}^d(\sigma_i-1)q_i
\end{equation}
in $R(K)^\omega$.
We put $\lambda_1=N_{K/k(\mu_p)}\Lambda_1=N_G\Lambda_1$.
Then $\lambda_1\mathcal{O}_{k(\mu_p)}$ is a degree one prime
ideal of $k(\mu_p)$.
Let $\lambda_1^\omega\in\mathcal{O}_{k(\mu_p),p}^\times$
be an element such that the projection of 
$\lambda_1\otimes 1\in\mathcal{O}_{K,p}^\times\otimes\F_p$
to the $\omega$-part is $\lambda_1^\omega\otimes 1
\in(\mathcal{O}_{K,p}^\times\otimes\F_p)^\omega$.
Then it follows from
\eqref{c1} and \eqref{c2} that
$
\varepsilon_0\lambda_1^\omega\otimes 1
\in(\mathcal{O}_{K,p}^\times\otimes\F_p)^\omega
$
maps to 0 by the natural map 
$
(\mathcal{O}_{K,p}^\times\otimes\F_p)^\omega
\longrightarrow
R(K)^\omega
$
because
$(\lambda_1^\omega \mod U(K)^p)
=(\lambda_1 \mod U(K)^p)=N_G(\Lambda_1 \mod U(K)^p)
=-(\varepsilon_0\mod U(K)^p)$ in $R(K)^\omega$.
Hence every prime of $L_p(k)$
lying over $p$ splits completely and
that dividing $N_{k(\mu_p)/k}\lambda_1$ totally ramifies
in the degree $p$ extension
$L_p(k)\{\sqrt[p]{\varepsilon_0\lambda_1^\omega}\}/
L_p(k)$, respectively.
By using Lemma \ref{lem-9} and \eqref{td} again, we can choose
a degree one prime element $\Lambda_2\in\mathcal{O}_{K}$ which is prime to $p$
such that
\begin{equation}\label{c3}
(\Lambda_2\mod U(K)^p)=-x+\sum_{i=1}^d(\sigma_i-1)q_{i+d}
\end{equation}
in $R(K)^\omega$ and that the prime $\Lambda_2\mathcal{O}_{K}$ of
$K$ remains prime
in $K(\sqrt[p]{\varepsilon_0\lambda_1^\omega})$.
Let $\lambda_2=N_{K/k(\mu_p)}\Lambda_2$.
Then, similarly to the above, we see that $\lambda_2$
is also a prime element of $\mathcal{O}_{k(\mu_p)}$, and
that every prime of $L_p(k)$
lying over $p$ splits completely and
that dividing $N_{k(\mu_p)/k}\lambda_2$ totally ramifies
in the degree $p$ extension $L_p(k)\{\sqrt[p]{\varepsilon_0\lambda_2^\omega}\}/
L_p(k)$,
$\lambda_2^\omega\in\mathcal{O}_{K,p}^\times$
being similar to $\lambda_1^\omega$, respectively.
\par
In what follows, we shall show that 
$L_p(k)\{\sqrt[p]{\varepsilon_0\lambda_1^\omega},\ 
\sqrt[p]{\varepsilon_0\lambda_2^\omega}\}$
is the maximal elementary abelian $p$-extension of $L_p(k)$ which
is unramified outside the primes dividing
$N_{k(\mu_p)/k}(\lambda_1\lambda_2)$.

Assume that $L_p(k)\{\sqrt[p]{\beta}\}/L_p(k)$ is 
unramified outside the primes dividing
$N_{k(\mu_p)/k}(\lambda_1\lambda_2)$ for 
$\beta\otimes 1\in (\mathcal{O}_{K,p}^\times\otimes\F_p)^\omega$.
Then we find that 
\begin{equation}\label{beta}
\beta\otimes 1
=
\eta\otimes 1+\sum_{\sigma\in G}a_\sigma\sigma(\Lambda_1^\omega\otimes 1)
+\sum_{\sigma\in G}b_\sigma\sigma(\Lambda_2^\omega\otimes 1)
\ \ \mbox{in}\ \ (\mathcal{O}_{K,p}^\times\otimes\F_p)^\omega,
\end{equation}
where $\Lambda_i^\omega\otimes 1$
with $\Lambda_i^\omega\in\mathcal{O}_{K,p}^\times$
is the projection of $\Lambda_i\otimes 1\in\mathcal{O}_{K,p}^\times\otimes\F_p$ to
$(\mathcal{O}_{K,p}^\times\otimes\F_p)^\omega\ (i=1,2)$,
for some $\eta\otimes 1\in(\mathcal{O}_{K}^\times\otimes\F_p)^\omega$
and $a_\sigma, b_\sigma\in\F_p$.
Since all the primes of $K$ lying over $p$ are unramified in 
$K(\sqrt[p]{\beta})/K$,
we have
\begin{equation*}
\begin{aligned}
(\eta\mod U'(K))&+\sum_{\sigma\in G}a_\sigma\sigma(\Lambda_1^\omega
\mod U'(K))\\
&+\sum_{\sigma\in G}b_\sigma\sigma(\Lambda_2^\omega
\mod U'(K))
=0
\end{aligned}
\end{equation*}
in $R'(K)^\omega$,
which implies
\begin{align*}
(\eta&\mod U(K)^p)-\sum_{\sigma\in G}(a_\sigma+b_\sigma)\sigma x\\
&+\sum_{\sigma\in G}a_\sigma\sigma\sum_{i=1}^d(\sigma_i-1)q_{i}
+\sum_{\sigma\in G}b_\sigma\sigma\sum_{i=1}^d(\sigma_i-1)q_{d+i}
\in Z
\end{align*}
in $R(K)^\omega$ by \eqref{c2}, \eqref{c3}, and the fact
$(\Lambda_i\mod U(K)^p)=(\Lambda_i^\omega\mod U(K)^p)$
in $R(K)^\omega$ for $i=1,2$.
Hence we obtain
\begin{equation}\label{d2}
(\eta\mod U(K)^p)-\sum_{\sigma\in G}(a_\sigma+b_\sigma)\sigma x
=0,
\end{equation}
and
\begin{equation}\label{d1}
\sum_{\sigma\in G}a_\sigma\sigma(\sigma_i-1)
=\sum_{\sigma\in G}b_\sigma\sigma(\sigma_{i}-1)
=0\ \ (1\le i\le d),
\end{equation}
because
$(\eta\mod U(K)^p)-\sum_{\sigma\in G}(a_\sigma+b_\sigma)\sigma x
\in P$ and $q_i$'s form a system of free basis of $Q$ over $\F_p[G]$.
It follows from \eqref{d1} that
there exist some
$a, b\in\F_p$ such that
$a_\sigma=a$ and
$b_\sigma=b$ for all $\sigma\in G$,
because $\sum_{\sigma\in G}a_\sigma\sigma(\tau-1)
=\sum_{\sigma\in G}b_\sigma\sigma(\tau-1)=0$
holds for any $\tau\in G$ by \eqref{d1}
and our choice of $\sigma_i$'s,
which implies
$\sum_{\sigma\in G}a_\sigma\sigma,\ 
\sum_{\sigma\in G}b_\sigma\sigma\in\F_p[G]^G=\F_pN_G$.
Hence we derive from \eqref{c1} and \eqref{d2} that
\[
(\eta\mod U(K)^p)-(a+b)(\varepsilon_0\mod U(K)^p)=0,
\]
which implies
$\eta\otimes 1=\varepsilon_0^{a+b}\otimes 1$
in $(\mathcal{O}_K^\times\otimes\F_p)^\omega$
by the injectivity of 
$(\mathcal{O}_K^\times\otimes\F_p)^\omega\longrightarrow R(K)^\omega.$
Therefore, by using \eqref{beta} and the fact 
$N_G(\Lambda_i^\omega\otimes 1)=\lambda_i^\omega\otimes 1\ (i=1,2)$,
we conclude that
\[
\beta\otimes 1=\varepsilon_0^{a+b}(\lambda_1^\omega)^a(\lambda_2^\omega)^b
\otimes 1
\ \ \mbox{in}\ \ (\mathcal{O}_{K,p}^\times\otimes\F_p)^\omega,
\]
and that
\[
L_p(k)\{\sqrt[p]{\beta}\}=L_p(k)\{\sqrt[p]{
(\varepsilon_0\lambda_1^\omega)^a(\varepsilon_0\lambda_2^\omega)^b}\}\subseteq
L_p(k)\{\sqrt[p]{\varepsilon_0\lambda_1^\omega},
\ \sqrt[p]{\varepsilon_0\lambda_2^\omega}\}.
\]
Thus we have shown that 
$L_p(k)\{\sqrt[p]{\varepsilon_0\lambda_1^\omega},
\ \sqrt[p]{\varepsilon_0\lambda_2^\omega}\}$
is the maximal elementary abelian $p$-extension of $L_p(k)$
which is unramified outside the primes dividing
$N_{k(\mu_p)/k}(\lambda_1\lambda_2)$.
\par
Since the prime $\Lambda_2\mathcal{O}_{K}$ of $K$
remains prime in
$K(\sqrt[p]{\varepsilon_0\lambda_1^\omega})$
and ramifies in $K(\sqrt[p]{\varepsilon_0\lambda_2^\omega})$
by our choice of $\Lambda_1$ and $\Lambda_2$,
the decomposition subgroup of
$\Gal(L_p(k)\{\sqrt[p]{\varepsilon_0\lambda_1^\omega},
\ \sqrt[p]{\varepsilon_0\lambda_2^\omega}\}/L_p(k))$ at the prime
below $\Lambda_2\mathcal{O}_K$
is the whole of the Galois group.
Hence it follows from Lemma \ref{cent} that
$L_p(k)\{\sqrt[p]{\varepsilon_0\lambda_1^\omega},
\ \sqrt[p]{\varepsilon_0\lambda_2^\omega}\}$
has no non-trivial unramified $p$-extensions:
\begin{equation}\label{K}
L_p(L_p(k)\{\sqrt[p]{\varepsilon_0\lambda_1^\omega},
\ \sqrt[p]{\varepsilon_0\lambda_2^\omega}\})=
L_p(k)\{\sqrt[p]{\varepsilon_0\lambda_1^\omega},
\ \sqrt[p]{\varepsilon_0\lambda_2^\omega}\}.
\end{equation}
Also, every prime of $L_p(k)$ lying over $p$ splits completely
in the extension $L_p(k)\{\sqrt[p]{\varepsilon_0\lambda_1^\omega},
\ \sqrt[p]{\varepsilon_0\lambda_2^\omega}\}/L_p(k)$.
\par
Now, 
$k':=k\{\sqrt[p]{\lambda_1^\omega(\lambda_2^\omega)^{-1}}\}$
is a required field as we shall see in the following:
We first note that
$L_p(k)\{\sqrt[p]{\varepsilon_0\lambda_1^\omega}\}/k$
is also a solution of embedding problem \eqref{eprob}
by Theorem A,
and that every prime of $k'$ lying over $p$ splits completely in
$L_p(k)\{\sqrt[p]{\varepsilon_0\lambda_1^\omega}\},
\ \sqrt[p]{\varepsilon_0\lambda_2^\omega}\}/k'$.
Since $k'/k$ is linearly disjoint from
$L_p(k)\{\sqrt[p]{\varepsilon_0\lambda_1^\omega}\}/k$,
the restriction map induces the isomorphism
$\Gal(L_p(k)\{\sqrt[p]{\varepsilon_0\lambda_1},\ 
\sqrt[p]{\varepsilon_0\lambda_2}\}/k')\simeq
\Gal(L_p(k)\{\sqrt[p]{\varepsilon_0\lambda_1^\omega}\}/k)\simeq G'$.
Also, because 
$L_p(k)\{\sqrt[p]{\varepsilon_0\lambda_1^\omega},\ 
\sqrt[p]{\varepsilon_0\lambda_2^\omega}\}
/L_p(k)\{\sqrt[p]{\lambda_1^\omega(\lambda_2^\omega)^{-1}}\}$
and $L_p(k)\{\sqrt[p]{\lambda_1^\omega(\lambda_2^\omega)^{-1}})/k'$
are unramified $p$-extensions,
we conclude that $L_p(k)\{\sqrt[p]{\varepsilon_0\lambda_1^\omega},\ 
\sqrt[p]{\varepsilon_0\lambda_2^\omega}\}/k'$ is an unramified $p$-extension.
Hence it follows from \eqref{K}
that $L_p(k)\{\sqrt[p]{\varepsilon_0\lambda_1^\omega},\
\sqrt[p]{\varepsilon_0\lambda_2^\omega}\}=L_p(k')$.
Therefore we find that
$\t{G}_{k'}(p)$ fits commutative diagram \eqref{cd-propII}
with suitable isomorphism $\theta$, and that $L_p(k')/k'$ is $p$-decomposed.
\par
Finally, we shall examine the condition on $B_p(k')$
by using the following:
\begin{lem}(Iwasawa, \cite[p.7, (4)]{Iw73})\label{iwlem}
Let $F/E$ be a cyclic extension of prime degree $p$,
and let $s$ be the number ramified primes in $F/E$.
Then we have
\[
\dim_{\F_p}(\Cl(F)\otimes\F_p)\le
p(\dim_{\F_p}(\Cl(E)\otimes\F_p)+s).
\]
\end{lem}
By applying the above lemma to cyclic extensions
$L_p(k')(\mu_p)/Kk'$ and $Kk'/K$ of degree $p$,
we derive from
our assumption \eqref{B_p(k)} on $B_p(k)$
that
$r_2(k')\ge B_p(k')$;
Since the number of ramified primes at $Kk'/K$
is $2[K:k]$, which does not exceed $2(p-1)\#G$,
and $L_p(k')(\mu_p)/Kk'$ is unramified,
we have $\dim_{\F_p}(\Cl(L_p(k')(\mu_p))\otimes\F_p)
\le p^2(\dim_{\F_p}(\Cl(K)\otimes\F_p)+2(p-1))$
by Lemma \ref{iwlem}.
\par
Thus we obtain a required field $k'$. This proves Proposition II,
hence we have completed the proof of Theorems 1 and 2.

\par\bigskip\noindent
\hspace{6cm}Manabu Ozaki,\par\noindent
 \hspace{6cm}Department of Mathematics,\par\noindent
\hspace{6cm}School of Science and Engineering,\par\noindent
\hspace{6cm}Kinki University,\par\noindent
\hspace{6cm}Kowakae 3-4-1,\par\noindent
\hspace{6cm}Higashi-Osaka 577-8502, Japan\par\noindent
\hspace{6cm}e-mail:\ \verb+ozaki@math.kindai.ac.jp+
\end{document}